\documentclass{article}

\usepackage[utf8]{inputenc}
\usepackage{amsmath}
\usepackage{amssymb}
\usepackage{amsthm}
\usepackage{hyperref}
\usepackage{algorithm, algorithmic}
\usepackage{listings}
\usepackage{graphicx, caption, subcaption}
\graphicspath{{figures/}}
\usepackage{xcolor}
\usepackage{enumitem}
\usepackage{fullpage}
\usepackage{booktabs}
\usepackage{algorithm}
\usepackage{algorithmic}
\usepackage{natbib}

\newtheorem{theorem}{\bf Theorem}[section]
\newtheorem{definition}[theorem]{\bf Definition}
\newtheorem{lemma}[theorem]{\bf Lemma}
\newtheorem{remark}[theorem]{\bf Remark}

\newtheorem{corollary}[theorem]{\bf Corollary}

\newtheorem{assumption}{Assumption}

\newcommand{\R}{\mathbb{R}}

\newcommand{\Rmm}{\mathbb{R}^{m \times m}}

\newcommand{\norm}[1]{\left\lVert#1\right\rVert}
\newcommand{\normStd}[1]{\lVert#1\rVert}
\newcommand{\Mr}{\mathcal{M}_r}

\newcommand{\Tq}{\mathcal{T}_q}
\newcommand{\Tr}{\mathcal{T}_r}

\newcommand{\Trp}{\mathcal{T}_r^{\perp}}
\newcommand{\flow}{\phi^h}
\newcommand{\flowq}{\psi^h_q}
\newcommand{\flowr}{\psi^h_r}
\newcommand{\Xnk}{X_n^k}
\newcommand{\Enk}{E_n^k}
\newcommand{\Xnkk}{X_n^{k+1}}
\newcommand{\Enkk}{E_n^{k+1}}
\newcommand{\Xnnkk}{X_{n+1}^{k+1}}
\newcommand{\Ennkk}{E_{n+1}^{k+1}}
\newcommand{\suminf}{\sum_{n=0}^{\infty}}
\newcommand{\suminfone}{\sum_{n=1}^{\infty}}
\newcommand{\proj}[2]{\mathcal{P}_{#1}#2}
\newcommand{\scal}[1]{\left\langle #1 \right\rangle}

\DeclareMathOperator{\rank}{rank}
\DeclareMathOperator{\diag}{diag}

\usepackage{accents}
\newcommand*{\dt}[1]{%
  \accentset{\mbox{\large\bfseries .}}{#1}}

\title{Low-rank Parareal:\\ a low-rank parallel-in-time integrator}

\author{Benjamin Carrel, Martin J.~Gander, Bart Vandereycken\thanks{Section of Mathematics, University of Geneva. This work was supported by the SNSF under research project 192363.}}
\date{September 5, 2022}

\begin{document}

\maketitle

\begin{abstract}
  In this work, the Parareal algorithm is applied to evolution problems that admit good low-rank approximations and for which the dynamical low-rank approximation (DLRA) can be used as time stepper. Many discrete integrators for DLRA have recently been proposed, based on splitting the projected vector field or by applying projected Runge--Kutta methods. The cost and accuracy of these methods are mostly governed by the rank chosen for the approximation. These properties are used in a new method, called low-rank Parareal, in order to obtain a time-parallel DLRA solver for evolution problems. The algorithm is analyzed on affine linear problems and the results are illustrated numerically.
\end{abstract}

\section{Introduction} \label{sec:introduction}

This work is concerned with the parallel-in-time integration of evolution problems for which the solution can be well approximated by a time-dependent low-rank matrix. In particular, we aim to solve approximately the evolution problem
\begin{align} \label{general_problem}
  \begin{aligned}
    \dt{X}(t) & = F ( t, X(t) ), &  & t \in [0,T], \\
    X(0)      & = X_0,
  \end{aligned}
\end{align}
where $X(t)$ is a matrix of size $m \times m$. When the dimension $m$ is large, the numerical solution of~\eqref{general_problem} can be very expensive since the matrix $X(t)$ is usually dense. One way to alleviate this curse of dimensionality is to use low-rank approximations where, for every $t$, we approximate $X(t)$ by $Y(t) \in \Rmm$ such that $\rank(Y(t)) = r \ll m$. The accuracy of this approximation will depend on the choice of the rank $r$. Here, $X(t)$ is assumed to be square for notational convenience and all results can be easily formulated for rectangular $X(t)$.

A popular paradigm to solve directly for the low-rank approximation $Y(t)$ is the dynamical low-rank algorithm (DLRA), first proposed in~\cite{koch_dynamical_2007}. As defined later in Def.~\ref{DLRA}, DLRA leads to an evolution problem that is a projected version of \eqref{general_problem}. In the last decade, many discrete integrators for this projected problem have been proposed. One class of integrators consists in a clever splitting of the projector so that the resulting splitting method can be implemented efficiently. An influential example is the projector-splitting scheme proposed in~\cite{lubich_projector-splitting_2014}. Other methods that require integrating parts of the vector field can be found in~\cite{khoromskij_efficient_2012,cerutiUnconventionalRobustIntegrator2022}. Another approach, proposed in~\cite{feppon_dynamically_2018,kieri_projection_2019,rodgers_adaptive_2021}, is based on projecting standard Runge--Kutta methods (sometimes including their intermediate stages). Most of these methods are formulated for constant rank $r$.  Rank adaptivity can be incorporated without much difficulty for splitting  and for projected schemes; see~\cite{feppon_dynamically_2018,dektor_rank-adaptive_2021,rodgers_adaptive_2021,cerutiRankadaptiveRobustIntegrator2022}. Finally, given the importance of DLRA in problems from physics (like the Schrödinger and Vlasov equation), the integrators in~\cite{lubich_projector-splitting_2014,einkemmerQuasiConservativeDynamicalLowRank2019} also preserve certain invariants, like energy. However, none of these time integrators consider a parallel-in-time scheme for DLRA, which is particularly interesting in the large-scale setting.




Parallel computing can be very effective and is even necessary to solve (very) large-scale problems. While parallelization in space is well known, also the time direction can be parallelized to some extent when solving evolution problems.  Over the last decades, various parallel-in-time algorithms have been proposed; see, e.g., the overviews~\cite{gander50YearsTime2015a,ongApplicationsTimeParallelization2020}. Among these, the Parareal algorithm from~\cite{lions_resolution_2001} is one of the more popular algorithms for time parallelization. It is based on a Newton-like iteration, with inaccurate but cheap corrections performed sequentially, and accurate but expensive solves performed in parallel. This idea of solving in parallel an evolution problem as a nonlinear (discretized) system also appears in related methods like PFASST~\cite{emmettEfficientParallelTime2012}, MGRIT~\cite{friedhoffMultigridintimeAlgorithmSolving2012}  and Space-Time Multi-Grid~\cite{ganderAnalysisNewSpacetime2016}. Theoretical results and numerical studies on a large numbers of cores show that these parallel-in-time methods can have good parallel performance for parabolic problems; see, e.g., \cite{speckMassivelySpacetimeParallel2012,hoferParallelRobustPreconditioning2019,ganderUnifiedAnalysisFramework2022}. So far, these methods did not incorporate a low-rank compression of the space dimension, which is the main topic of this work.

\section{Preliminaries and contributions}

\subsection{The Parareal algorithm} \label{sec:parareal_algo}

The Parareal iteration in Def.~\ref{def:parareal} below is given for constant time step $h$ (hence $T=Nh$) and for autonomous $F$. Both restrictions are not crucial but ease the presentation. The quantity $X_n^k$ is an approximation for $X(t_n)$ at time $t_n = nh$ and iteration $k$. The accuracy of this approximation is expected to improve with increasing $k$. Here, and throughout the paper, we denote dependency on the iteration index $k$ as ${\ }^k$, which should not be confused with the $k$th power.

\begin{definition}[Parareal]\label{def:parareal}
  The Parareal algorithm is defined by the following double iteration on $k$ and $n$,
  \begin{align} \label{def:Parareal}
     & \text{(Initial value)} \quad          &  & X_0^k     = X_0,                                                              \\
     & \text{(Initial approximation)}  \quad &  & X_{n+1}^0 = \mathcal G^h (X_n^0),                                             \\
     & \text{(Iteration)} \quad              &  & \Xnnkk    = \mathcal F^h (\Xnk) + \mathcal G^h (\Xnkk) - \mathcal G^h (\Xnk).
  \end{align}
  Here, $\mathcal F^h(X)$ represents a fine (accurate) time stepper applied to the initial value $X$ and propagated until time $h$. Similarly, $\mathcal G^h(X)$ represents a coarse (inaccurate) time stepper.
\end{definition}

Given two time steppers, the Parareal algorithm is easy to implement.
A remarkable property of Parareal is the convergence in a finite number of steps for $k=n$. It is well known that Parareal works well on parabolic problems but behaves worse on hyperbolic problems; see~\cite{gander_analysis_2007} for an analysis.

\subsection{Dynamical low-rank approximation} \label{sec:dlra}

Let $\Mr$ denote the set of $m \times m$ matrices of rank $r$, which is a smooth embedded submanifold in $\Rmm$. Instead of solving~\eqref{general_problem}, the DLRA solves the following projected problem:

\begin{definition}[Dynamical low-rank approximation]
  For a rank $r$, the dynamical low-rank approximation of problem~\eqref{general_problem} is the solution of
  \begin{align} \label{DLRA}
    \begin{aligned}
       & \dt{Y}(t) = \proj{Y(t)}{F(t, Y(t))}, \\
       & Y(0) = Y_0 \in \Mr,
    \end{aligned}
  \end{align}
  where $\mathcal{P}_{Y}$ is the $l_2$-orthogonal projection onto the tangent space $\mathcal T_{Y} \Mr$ of $\Mr$ at $Y \in \Mr$. In particular, $Y(t) \in \Mr$ for every $t \in [0,T]$.
\end{definition}

To analyze the approximation error of DLRA, we need the following standard assumptions from~\cite{kieri_discretized_nodate}. Here, and throughout the paper, $\| \cdot \|$ denotes the Frobenius norm.
\begin{assumption}[DLRA assumptions]\label{ass:DLRA assumptions}
  The function $F$ satisfies the following properties for all $X,Y \in \Rmm$:
  \begin{itemize}
    \item Lipschitz with constant $L$: $\norm{F(X)-F(Y)} \leq L \norm{X-Y}$.
    \item One-sided Lipschitz with constant $\ell$: $\langle X-Y,F(X)-F(Y) \rangle \leq \ell \norm{X-Y}^2$.
    \item Maps almost to the tangent bundle of $\Mr$: $\norm{F(Y)-\mathcal{P}_Y F(Y)} \leq \varepsilon_r$.
  \end{itemize}
\end{assumption}

In the analysis in Section~\ref{sec:low-rank_Parareal}, it is necessary to have $\ell < 0$ for convergence. This holds when $F$ is a discretization of certain parabolic PDEs, like the heat equation. In particular, for an affine function of the form $F(X)=A(X)+B$, it holds
$\ell = \tfrac{1}{2}\lambda_{\max}(A+A^T)$; see~\cite[Ch.~I.10]{hairer_solving_1987}. The quantity $\varepsilon_r$ is called the \textit{modeling error} and decreases when the rank $r$ increases. For our problems of interest, this quantity is typically very small. Finally, the existence of $L$ is only needed to guarantee the uniqueness of~\eqref{general_problem} but it will actually not appear in our analysis. We can therefore allow $L$ to be large, as is the case for discretized parabolic PDEs.

Standard theory for perturbations of ODEs allows us to obtain the following error bound from the assumptions above:

\begin{theorem}[Error of DLRA~\cite{kieri_discretized_nodate}]\label{thm:error of DLRA}
  Under Assumption~\ref{ass:DLRA assumptions}, the DLRA verifies
  \begin{align}
    \norm{\flowr(Y_0) - \flow(X_0)} \leq e^{\ell t} \norm{Y_0 - X_0} + \varepsilon_r \int_0^t e^{\ell s} ds,
  \end{align}
  where $\flow$ is the flow of the original problem~\eqref{general_problem} and $\flowr$ is the flow of its DLRA~\eqref{DLRA} for rank $r$.
\end{theorem}

The solution of DLRA~\eqref{DLRA} is quasi-optimal with the best rank approximation. This can be seen already in Theorem~\ref{thm:error of DLRA} for short time intervals. Similar estimates exist for parabolic problems~\cite{conte_dynamical_2020} and for longer time when there is a sufficiently large gap in the singular values and when their derivatives are bounded~\cite{koch_dynamical_2007}.

\subsection{Contributions} \label{sec:contributions}

In this paper, we propose a new algorithm, called \emph{low-rank Parareal}. As far as we know, this is the first parallel-in-time integrator for low-rank approximations. We analyze the proposed algorithm when the function $F$ in~\eqref{general_problem} is affine. To this end, we extend the analysis of the classical Parareal algorithm in \cite{barth_nonlinear_2008} to a more general setup where the coarse problem is different from the fine problem. We can prove that the method converges for big steps (large $h$) on diffusive problems ($\ell<0$). In numerical experiments, we confirm this behavior. In addition, the method also performs well empirically with a less strict condition on $h$ and on a non-affine problem.

\section{Low-rank Parareal} \label{sec:low-rank_Parareal}

We now present our low-rank Parareal algorithm for solving~\eqref{general_problem}. Since the cost of most discrete integrators for DLRA scales quadratically\footnote{While the actual cost can be larger, at the very least the methods typically compute compact QR factorizations to normalize the approximations. This costs $O(mr^2 + r^3)$ flops in our setting.} with the approximation rank, we take the coarse time stepper as DLRA with a small rank $q$. Likewise, the fine time stepper is DLRA with a large rank $r$. We can even take $r=m$, which corresponds to computing the exact solution as the fine time stepper since $Y \in \mathbb{R}^{m \times m}$.

\begin{definition}[Low-rank Parareal] \label{def:low_rank_parareal}
  Consider two ranks $q < r$. The low-rank Parareal algorithm iterates
  \begin{align}
     & \text{(Initial value)} \quad          &  & Y_0^k = Y_0, \label{def:initial value}                                                                                      \\
     & \text{(Initial approximation)}  \quad &  & Y_{n+1}^0 = \flowq \circ \Tq (Y_n^0) + \mathcal E_n, \label{def:initial approximation}                                                              \\
     & \text{(Iteration)} \quad              &  & Y_{n+1}^{k+1} = \flowr \circ \Tr(Y_n^k) + \flowq \circ \Tq (Y_n^{k+1}) - \flowq \circ \Tq (Y_n^k), \label{def:iteration}
  \end{align}
  where $\flowr(Z)$ is the solution of~\eqref{DLRA} at time $h$ with initial value $Y_0 = Z$, and $\Tr$ is the orthogonal projection onto $\Mr$. The notations $\flowq$ and $\Tq$ are similar but apply to rank $q$. The matrices $\mathcal E_n$ are small perturbations such that $\rank(Y_{n+1}^0) = r + 2q$ and can be chosen randomly.\footnote{One could even take all initial $Y_{n+1}^0$ random, as is sometimes also done in standard Parareal. The important property is $\rank(Y_{n+1}^0) = r + 2q$ so that the low-rank Parareal iterations~\eqref{def:iteration} are performed on the correct manifold.} 
\end{definition}

Observe that the rank of $Y_n^k$ is at most $r+2q$ for all $n,k$. The low-rank structure is therefore preserved over the iterations.
The matrices $\mathcal E_n$ insure that each iteration has a rank between $r$ and $r+2q$. These matrices impact only the initial error but do not have any role in the convergence of the algorithm as is shown later in the analysis.
An efficient implementation should store the low-rank matrices in factored form. In this context, the truncated SVD can be efficiently performed. The DLRA flows $\flowr$ and $\flowq$ can only be computed for relatively small problems. For larger problems, a suitable DLRA integrator must be used; see Section \ref{sec:numerical_experiments} for implementation details.

\subsection{Convergence analysis} \label{sec:convergence_analysis}

Let $X_{n} = X(t_{n})$ be the solution of the full problem~\eqref{general_problem} at time $t_{n}$. Let $Y_{n}^{k}$ be the corresponding low-rank Parareal solution at iteration $k$. We are interested in bounding the error of the algorithm,
\begin{equation} \label{def:error}
  \begin{aligned}
    \Enk  = X_{n} - Y_{n}^{k},
  \end{aligned}
\end{equation}
for all relevant $n$ and $k$. To this end, we make the following assumption:
\begin{assumption}[Affine vector field]\label{ass:affine}
  The function $F$ is affine linear and autonomous, that is, $$F(X) = A(X) + B$$ with $A\colon \Rmm \to \Rmm$ a linear operator and $B \in \Rmm$.
\end{assumption}

The following lemma gives us a recursion for the Frobenius norm of the error. This recursion will be fundamental in deriving our convergence bounds later on when we generalize the proof for standard Parareal from~\cite{barth_nonlinear_2008}.

\begin{lemma}[Iteration of the error] \label{lem:iteration_error}
  Under the Assumptions~\ref{ass:DLRA assumptions} and~\ref{ass:affine}, the error of low-rank Parareal verifies
  \begin{align}\label{eq:error recursion detailed}
    \norm{\Ennkk} \leq e^{\ell h} C_{r,q} \norm{\Enk} + e^{\ell h} C_q \norm{\Enkk} + e^{\ell h} \max_{n \geq 0} \norm{X_n - \Tr(X_n)} + (2 \varepsilon_q + \varepsilon_r) \int_0^h e^{\ell (h-s)} ds.
  \end{align}
  The constants $\ell, \varepsilon_q, \varepsilon_r$ are defined in Assumption \ref{ass:DLRA assumptions}. Moreover, $C_{r,q}$ and $C_q$ are the Lipschitz constants of $\mathcal T_{r} - \mathcal T_{q}$ and $\Tq$.
\end{lemma}

\begin{proof}
  Our proof is similar to the one in \cite{kieri_projection_2019} where first the continuous version of the approximation error of DLRA is studied. Denote by $\phi^h(Z)$ the solution of~\eqref{general_problem} at time $h$ with initial value $Y_0 = Z$. By definition, the discrete error is
  \[
    \Ennkk = \phi^h(X_{n}) - \flowr \circ \Tr(Y_n^k) - \flowq \circ \Tq (Y_n^{k+1}) + \flowq \circ \Tq (Y_n^k).
  \]
  We can interpret each term above as a flow from $t_n$ to $t_{n+1}=t_n+h$. Denote these flows by $X(t), Z(t), W(t)$, and $V(t)$ with the initial values
  \begin{equation}\label{eq:intial value XZWV}
    X(t_n) = X_n, \ Z(t_n) = \Tr(Y_n^k), \ W(t_n) = \Tq (Y_n^{k+1}), \ V(t_n) = \Tq (Y_n^k).
  \end{equation}
  Defining the continuous error as
  $$ E(t) = X(t) - Z(t) - W(t) + V(t), $$
  we then get the identity $\Ennkk = E(t_{n}+h)$.

  We proceed by bounding $\norm{E(t)}$. By definition of the flows above, we have (omitting the dependence on $t$ in the notation)
  \begin{align*}
    \dt{E} & = F(X) - \proj{Z}{F(Z)} - \proj{W}{F(W)} + \proj{V}{F(V)}                                   \\
           & = F(X - Z - W + V) + F(Z) - \proj{Z}{F(Z)} - \proj{W}{F(W)} + F(W) - F(V) + \proj{V}{F(V)},
  \end{align*}
  where the last equality holds since the function $F$ is affine. Using Assumption~\ref{ass:DLRA assumptions} and Cauchy--Schwarz, we compute
  \begin{align*}
    \frac{1}{2} \frac{d}{dt} \norm{E(t)}^2 & = \scal{E, \dt{E}}                                                                                                     \\
                                           & = \scal{E, F(E)} + \scal{E, F(Z) - \proj{Z}{F(Z)}} + \scal{E, F(W) - \proj{W}{F(W)}} - \scal{E, F(V) - \proj{V}{F(V)}} \\
                                           & \leq \ell \norm{E}^2 + \varepsilon_r \norm{E} + 2 \varepsilon_q \norm{E}.
  \end{align*}
  Since $\frac{d}{dt} \norm{E(t)}^2 = 2 \norm{E(t)} \frac{d}{dt} \norm{E(t)}$, we therefore obtain the differential inequality
  \begin{align*}
    \frac{d}{dt} \norm{E(t)} \leq \ell \norm{E} + \varepsilon_r + 2 \varepsilon_q.
  \end{align*}
  Grönwall's lemma allows us to conclude
  \begin{align} \label{eq:iteration error inequality one}
    \norm{E(t_n+h)} \leq \norm{E(t_n)} e^{\ell h} + (2 \varepsilon_q + \varepsilon_r) \int_{t_n}^{t_n+h} e^{\ell (h-s)} ds.
  \end{align}
  From~\eqref{eq:intial value XZWV}, we get
  \begin{align*}
    E(t_n) & = X_n - \Tr(Y_n^k) - \Tq(Y_n^{k+1}) + \Tq(Y_n^k).
  \end{align*}
  Denoting $\Trp = I - \Tr$ and $\mathcal T_{r,q} = \Tr - \Tq$, we get after rearranging terms
  \begin{align*}
    E(t_n) & = \Trp(X_n) + \Tr(X_n) - \Tq(X_n) - \Tr(Y_n^k) + \Tq(Y_n^k) + \Tq(X_n) - \Tq(Y_n^{k+1})    \\
           & = \Trp(X_n) + \mathcal T_{r,q}(X_n) - \mathcal T_{r,q}(Y_n^k) + \Tq(X_n) - \Tq(Y_n^{k+1}).
  \end{align*}
  Taking norms gives
  \begin{equation} \label{eq:iteration error inequality two}
    \begin{aligned}
      \norm{E(t_n)} & \leq \norm{\Trp(X_n)} + \norm{\mathcal T_{r,q}(X_n) - \mathcal T_{r,q}(Y_n^k)} - \norm{\Tq(X_n) - \Tq(Y_n^{k+1})} \\
                    & \leq \max_{n \geq 0} \norm{\Trp(X_n)} + C_{r,q} \norm{\Enk} + C_q \norm{\Enkk},
    \end{aligned}
  \end{equation}
  where $C_{r,q}$ and $C_q$ are the Lipschitz constants of $\mathcal T_{r,q}$ and $\Tq$ respectively. Combining inequalities \eqref{eq:iteration error inequality one} and \eqref{eq:iteration error inequality two} gives the statement of the lemma.
\end{proof}

We now study the error recursion~\eqref{eq:error recursion detailed} in more detail. To this end, let us slightly rewrite it as
\begin{align}\label{eq:error recursion simplified}
  \norm{\Ennkk} \leq \alpha \norm{\Enk} + \beta \norm{\Enkk} + \kappa, \quad \norm{E_n^0} \leq \gamma,
\end{align}
with the non-negative constants
\begin{equation}\label{eq:def constants abck}
  \begin{aligned}
    \alpha & = e^{\ell h} C_{r,q}, \quad \beta = e^{\ell h} C_q, \quad \gamma = \max_{n \geq 0} \norm{E_n^0},                   \\
    \kappa & = e^{\ell h} \max_{n \geq 0} \norm{X_n - \Tr(X_n)} + (2 \varepsilon_q + \varepsilon_r) \int_0^h e^{\ell (h-s)} ds.
  \end{aligned}
\end{equation}

Our first result is a linear convergence bound, up to the DLRA approximation error. It is similar to the linear bound for standard Parareal.

\begin{theorem}[Linear convergence] \label{thm:linear_bound}
  Under the Assumptions~\ref{ass:DLRA assumptions} and~\ref{ass:affine}, and if $\alpha + \beta < 1$, low-rank Parareal verifies for all $k \in \mathbb{N}$ the linear bound
  \begin{align} \label{eq:linear_bound}
    \max_{n \geq 0} \norm{E_n^k} \leq \left( \frac{\alpha}{1 - \beta} \right)^k \max_{n\geq 0} \norm{E_n^0} + \frac{\kappa}{1 - \alpha - \beta},
  \end{align}
  where $\alpha, \beta, \kappa$ are defined in~\eqref{eq:def constants abck}.
\end{theorem}

\begin{proof}
  Define $e_{\star}^k = \max_{n \geq 0} \norm{\Enk}$. Taking the maximum for $n \geq 0$ of both sides of~\eqref{eq:error recursion simplified}, we obtain
  $$ e_{\star}^{k+1} \leq \alpha e_{\star}^k + \beta e_{\star}^{k+1} + \kappa. $$
  By assumption, $0 \le \beta < 1$ and we can therefore obtain the recursion
  $$ e_{\star}^{k+1} \leq \frac{\alpha}{1-\beta} e_{\star}^k + \frac{\kappa}{1-\beta}, $$
  with solution
  \[
    e_{\star}^k \leq \left(\frac{\alpha}{1 - \beta}\right)^k e_{\star}^0 + \frac{\kappa}{1- \alpha-\beta} \left[ 1 - \left(\frac{\alpha}{1  - \beta}\right)^k\right].
  \]
  By assumption, we also have $0 \leq \frac{\alpha}{1-\beta}<1$, which allows us to obtain the statement of the theorem.
\end{proof}

Next, we present a more refined superlinear bound. To this end, we require the following technical lemma that solves the equality version of the double iteration~\eqref{eq:error recursion simplified}. A similar result, but without the $\kappa$ term and only as an upper bound, already appeared in~\cite[Thm.~1]{barth_nonlinear_2008}. Our proof is therefore similar but more elaborate.

\begin{lemma} \label{lem:technical_lemma}
  Let $\alpha, \beta, \gamma, \kappa \in \R$ be any non-negative constants such that $\alpha < 1$ and $\beta < 1$. Let $e_n^k$ be a sequence depending on $n, k \in \mathbb{N}$ such that
  \begin{align}\label{eq:iteration technical lemma}
    e_{n+1}^{k+1}  = \alpha e_n^k + \beta e_n^{k+1} + \kappa, \qquad e_{n+1}^0      = \gamma, \qquad e_0^k = 0.
  \end{align}
  Then,
  \begin{equation}\label{eq:explicit solution enk}
    e_n^k = \kappa \sum_{j=0}^{k-1} \sum_{i=0}^{n-j-1} \binom{i+j}{i} \ \alpha^j \beta^i +
    \begin{cases}
      0                                                                & \text{ if } n \leq k,   \\
      \gamma \ \alpha^k  \sum_{i=0}^{n-k-1} \binom{i+k-1}{i} \ \beta^i & \text{ if } n \geq k+1.
    \end{cases}
  \end{equation}
\end{lemma}

\begin{proof}
  The idea is to use the generating function $\rho_k(\xi) = \suminfone e_n^k \xi^n$ for $k \geq 1$. Multiplying~\eqref{eq:iteration technical lemma} by $\xi^{n+1}$ and summing over $n$, we obtain
  \[
    \suminf e_{n+1}^{k+1} \xi^{n+1} = \suminf \alpha e_n^k \xi^{n+1} + \suminf \beta e_n^{k+1} \xi^{n+1} + \suminf \kappa \xi^{n+1}, \qquad \suminf e_{n+1}^0 \xi^{n+1} = \suminf \gamma \xi^{n+1}.
  \]
  Since $e_0^k = 0$ for all $k$, this gives the relations
  \[
    \rho_{k+1}(\xi) = \alpha \xi \rho_k (\xi) + \beta \xi \rho_{k+1} ( \xi) + \kappa  \frac{\xi}{1-\xi}, \qquad \rho_0(\xi) = \gamma \frac{\xi}{1 - \xi}.
  \]
  We can therefore obtain the linear recurrence
  \[
    \rho_{k+1}(\xi) =
    a \rho_k(\xi) + b, \quad \text{where } a = \frac{\alpha \xi}{1 - \beta \xi}, \ b = \frac{\kappa \xi}{(1-\xi)(1-\beta \xi)}.
  \]
  Its solution satisfies
  \[
    \rho_k(\xi) = \frac{\alpha^k \xi^k}{(1 - \beta \xi)^k} \frac{\gamma \xi}{1 - \xi} + \sum_{j=0}^{k-1} \frac{\alpha^j \xi^j}{(1-\beta \xi)^{j+1}} \frac{\kappa \xi}{1 - \xi}.
  \]

  It remains to compute the coefficients in the power series of the above formula since by definition of $\rho_k(\xi) = \suminfone e_n^k \xi^n$ they equal the unknowns $e_n^k$.
  The binomial series formula for $|z|<1$,
  \begin{equation}\label{eq:binomial series formula}
    \frac{1}{(1-z)^{k+1}} = \sum_{i=0}^{\infty} \binom{i + k}{i} z^i,
  \end{equation}
  together with the Cauchy product gives
  \begin{align*}
    \frac{1}{(1-\beta \xi)^k} \frac{1}{1-\xi}     & = \sum_{i=0}^{\infty} \binom{i+k-1}{i} \beta^i \xi^i \cdot \sum_{i=0}^{\infty} \xi^i = \suminf \left( \sum_{\ell=0}^n \binom{\ell+k-1}{\ell} \beta^{\ell} \right) \xi^n \\
    \frac{1}{(1-\beta \xi)^{j+1}} \frac{1}{1-\xi} & = \sum_{i=0}^{\infty} \binom{i+j}{i} \beta^i \xi^i \cdot \sum_{i=0}^{\infty} \xi^i = \suminf \left( \sum_{\ell=0}^n \binom{\ell+j}{\ell} \beta^{\ell} \right) \xi^n.
  \end{align*}
  Hence, the first term in $\rho_k(\xi)$ satisfies
  \begin{align*}
    \frac{\alpha^k \xi^k}{(1 - \beta \xi)^k} \frac{\gamma \xi}{1 - \xi}
     & = \gamma \alpha^k \suminf \left( \sum_{\ell=0}^n \binom{\ell+k-1}{\ell} \beta^{\ell} \right) \xi^{n+k+1}                 = \gamma \alpha^k \sum_{n=k+1}^{\infty} \left( \sum_{\ell=0}^{n-k-1} \binom{\ell+k-1}{\ell} \beta^{\ell} \right) \xi^n,
  \end{align*}
  while the second term can be written as
  \begin{align*}
    \sum_{j=0}^{k-1} \frac{\alpha^j \xi^j}{(1-\beta \xi)^{j+1}} \frac{\kappa \xi}{1 - \xi}
    = \kappa \sum_{j=0}^{k-1} \suminf \sum_{\ell=0}^{n} \binom{\ell+j}{\ell} \alpha^j \beta^{\ell} \xi^{n+j+1} = \kappa \sum_{n=1}^{\infty} \sum_{j=0}^{k-1} \left( \sum_{\ell=0}^{n-1} \binom{\ell+j}{\ell} \alpha^j \beta^{\ell} \right) \xi^{n+j}.
  \end{align*}
  Putting everything together, we have
  \begin{align*}
    \suminf e_n^k \xi^n = \gamma \alpha^k \sum_{n=k+1}^{\infty} \left( \sum_{\ell=0}^{n-k-1} \binom{\ell+k-1}{\ell} \beta^\ell \right) \xi^n + \kappa \sum_{m=1}^{\infty} \sum_{j=0}^{k-1} \left( \sum_{\ell=0}^{m-1} \binom{\ell+j}{\ell} \alpha^j \beta^{\ell} \right) \xi^{m+j}.
  \end{align*}
  Finally, we can identify the coefficient $e_n^k$ in front of $\xi^n$ with those on the right-hand side. The coefficient for $\xi^n$ in the first term is clearly nonzero only when $n \geq k+1$. In the second term, there is only one $m$ for every $j$ such that $m+j=n$. Substituting $m=n-j$ allows us to identify the coefficient of $\xi^n$.
\end{proof}

Using the previous lemma, we can obtain a convergence bound that is superlinear in $k$.

\begin{theorem}[Superlinear convergence] \label{thm:superlinear_bound}
  Under the Assumptions~\ref{ass:DLRA assumptions} and~\ref{ass:affine},
  and if $\alpha + \beta < 1$, the error of low-rank Parareal satisfies for all $n,k \in \mathbb{N}$ the bound
  \begin{align}
    \norm{E_n^k} \leq \frac{\alpha^k}{(k-1)!} \frac{\prod_{j=2}^k (n-j)}{1 - \beta} \max_{n \geq 0} \norm{E_n^0} + \frac{\kappa}{1 - \alpha - \beta},
  \end{align}
  where $\alpha, \beta, \kappa$ are defined in~\eqref{eq:def constants abck}.
\end{theorem}

\begin{proof}
  Define $e_n^k = \norm{\Enk}$. By Lemma~\ref{lem:iteration_error}, the terms $e_n^k$ verify the relation described in Lemma~\ref{lem:technical_lemma} with $=$ replaced by $\leq$ in~\eqref{eq:iteration technical lemma}. Hence, the solution~\eqref{eq:explicit solution enk} from Lemma~\ref{lem:technical_lemma} will be an upper bound for $e_n^k$.

  Since $0 \leq \alpha + \beta < 1$ and using the binomial series formula~\eqref{eq:binomial series formula}, we bound the first term in~\eqref{eq:explicit solution enk} as
  \begin{align*}
    \kappa \sum_{j=0}^{k-1} \sum_{i=0}^{n-j-1} \binom{i+j}{i} \alpha^j \beta^i
     & \leq \kappa  \sum_{j=0}^{k-1} \sum_{i=0}^{\infty} \binom{i+j}{i} \alpha^j \beta^i          = \kappa  \sum_{j=0}^{k-1} \alpha^j \frac{1}{(1-\beta)^{j+1}}                                  \\
     & \leq \frac{\kappa}{1-\beta} \sum_{j=0}^{\infty} \left( \frac{\alpha}{1-\beta} \right)^j  = \frac{\kappa}{1-\beta} \frac{1}{1-\frac{\alpha}{1-\beta}} = \frac{\kappa}{1 - \alpha - \beta}.
  \end{align*}
  For $0 \leq i \leq n-k-1$ and $n \geq k+1$, observe that
  \begin{align*}
    \frac{(i+k-1)!}{i!} = \prod_{j=1}^k (i+j) \leq \prod_{j=1}^k (n-k-1+j) = \prod_{j=2}^{k} (n-j).
  \end{align*}
  Since $0 \leq \beta < 1$, we can therefore bound the second term as
  \begin{align*}
    \gamma \ \alpha^k  \sum_{i=0}^{n-k-1} \binom{i+k-1}{i} \beta^i
     & = \gamma \ \alpha^k  \sum_{i=0}^{n-k-1} \frac{(i+k-1)!}{i! (k-1)!} \beta^i          \leq \gamma \ \frac{\alpha^k}{(k-1)!}  \prod_{j=2}^{k} (n-j) \sum_{i=0}^{n-k-1} \beta^i \\
     & \leq \gamma \ \frac{\alpha^k}{(k-1)!} \frac{\prod_{j=2}^k (n-j)}{1 - \beta}.
  \end{align*}
  The conclusion now follows by the definition of $\gamma$.
\end{proof}

The proof above can be modified to obtain a simple linear bound that is similar but different to the one from Theorem~\ref{thm:linear_bound}:

\begin{theorem}[Another linear convergence bound] \label{thm:second_linear_bound}
  Under Assumptions~\ref{ass:DLRA assumptions} and~\ref{ass:affine}, and if $\alpha + \beta < 1$, the error of low-rank Parareal satisfies for all $n,k \in \mathbb{N}$ the bound
  \begin{align} \label{eq:second_linear_bound}
    \norm{E_n^k} \leq \alpha^k (1+\beta)^{n-1} \max_{n \geq 0} \norm{E_n^0} + \frac{\kappa}{1 - \alpha - \beta},
  \end{align}
  where $\alpha, \beta, \kappa$ are defined in~\eqref{eq:def constants abck}.
\end{theorem}

\begin{proof}
  We repeat the proof for the superlinear bound but this time, the second term is bounded as
  \begin{equation*}
    \gamma \ \alpha^k  \sum_{i=0}^{n-k-1} \binom{i+k-1}{i} \beta^i \leq \gamma \ \alpha^k  \sum_{i=0}^{n-1} \binom{n-1}{i} \beta^i = \gamma \ \alpha^k \ (1+\beta)^{n-1}. \qed
  \end{equation*}
\end{proof}

\begin{remark}
  In the proof above, yet another bound based on~\eqref{eq:binomial series formula} is
  \begin{align*}
    \gamma \ \alpha^k  \sum_{i=0}^{n-k-1} \binom{i+k-1}{i} \beta^i \leq \gamma \ \alpha^k  \sum_{i=0}^{\infty} \binom{i+k-1}{i} \beta^i = \gamma \ \alpha^k \ \frac{1}{(1-\beta)^k}.
  \end{align*}
  This time we recover the linear bound from Theorem~\ref{thm:linear_bound}.
\end{remark}

\subsection{Summary of the convergence bounds}

\begin{figure}
  \centering
  \begin{subfigure}{0.49\textwidth}
    \includegraphics[width=\textwidth]{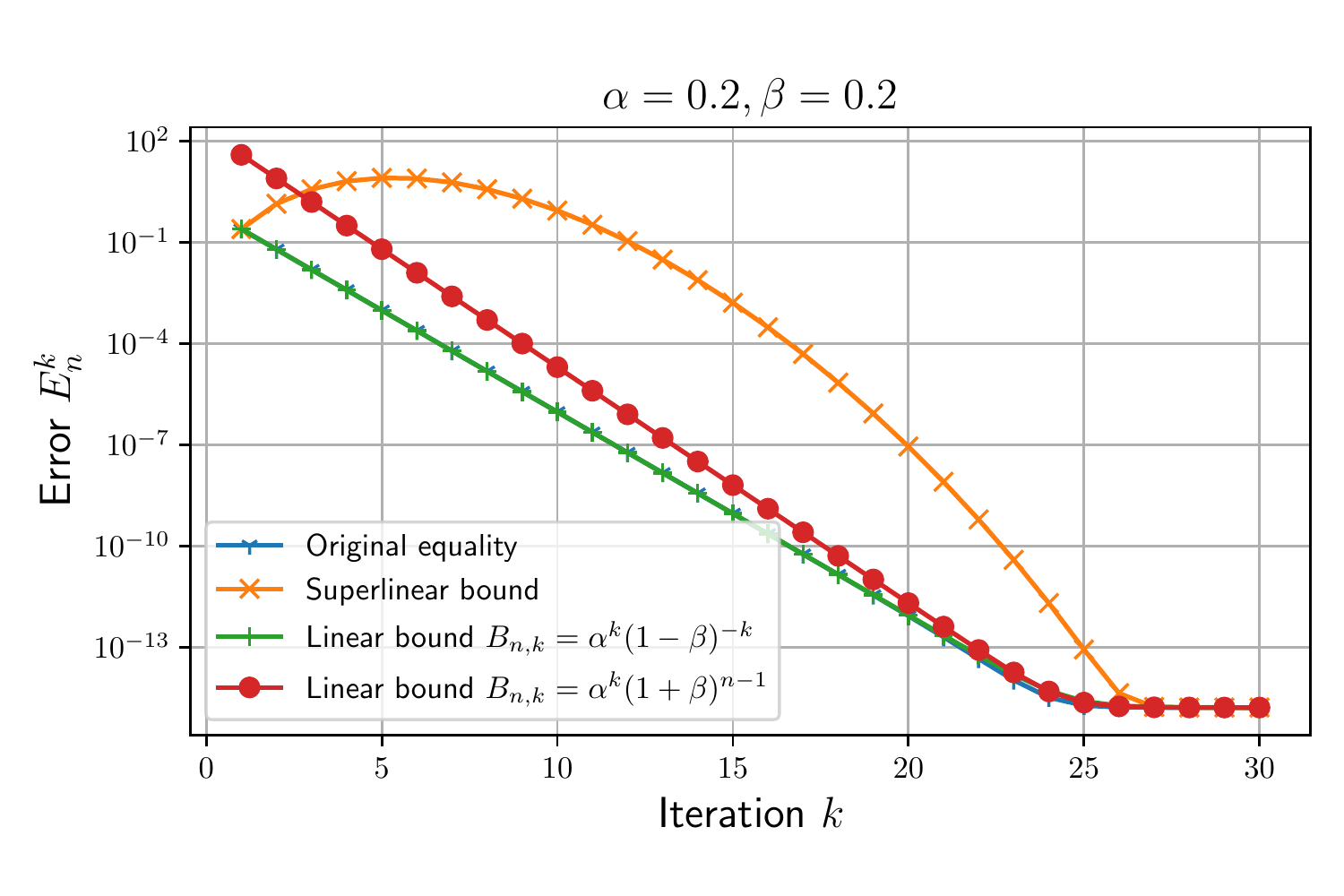}
  \end{subfigure}
  \begin{subfigure}{0.49\textwidth}
    \includegraphics[width=\textwidth]{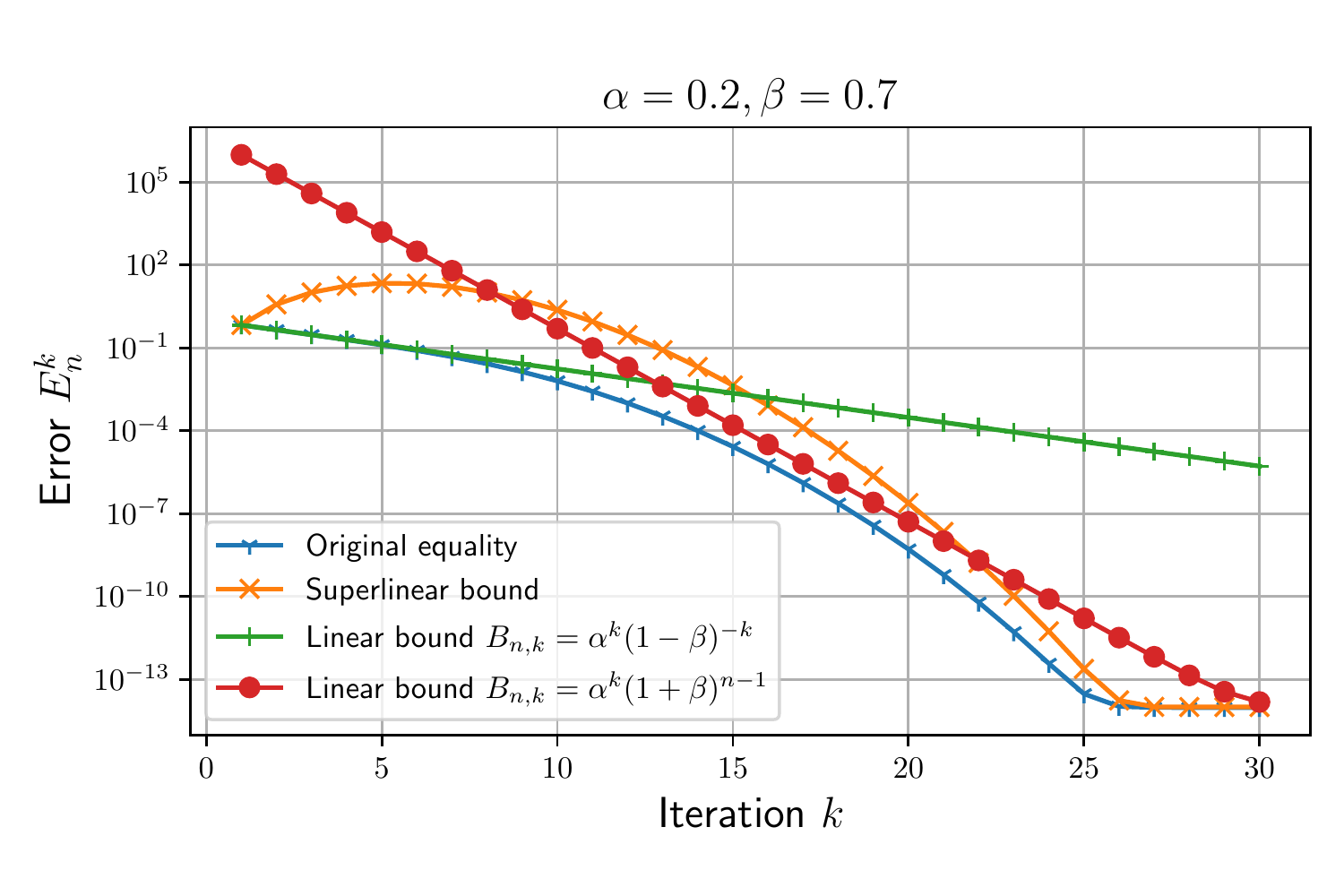}
  \end{subfigure}
  \begin{subfigure}{0.49\textwidth}
    \includegraphics[width=\textwidth]{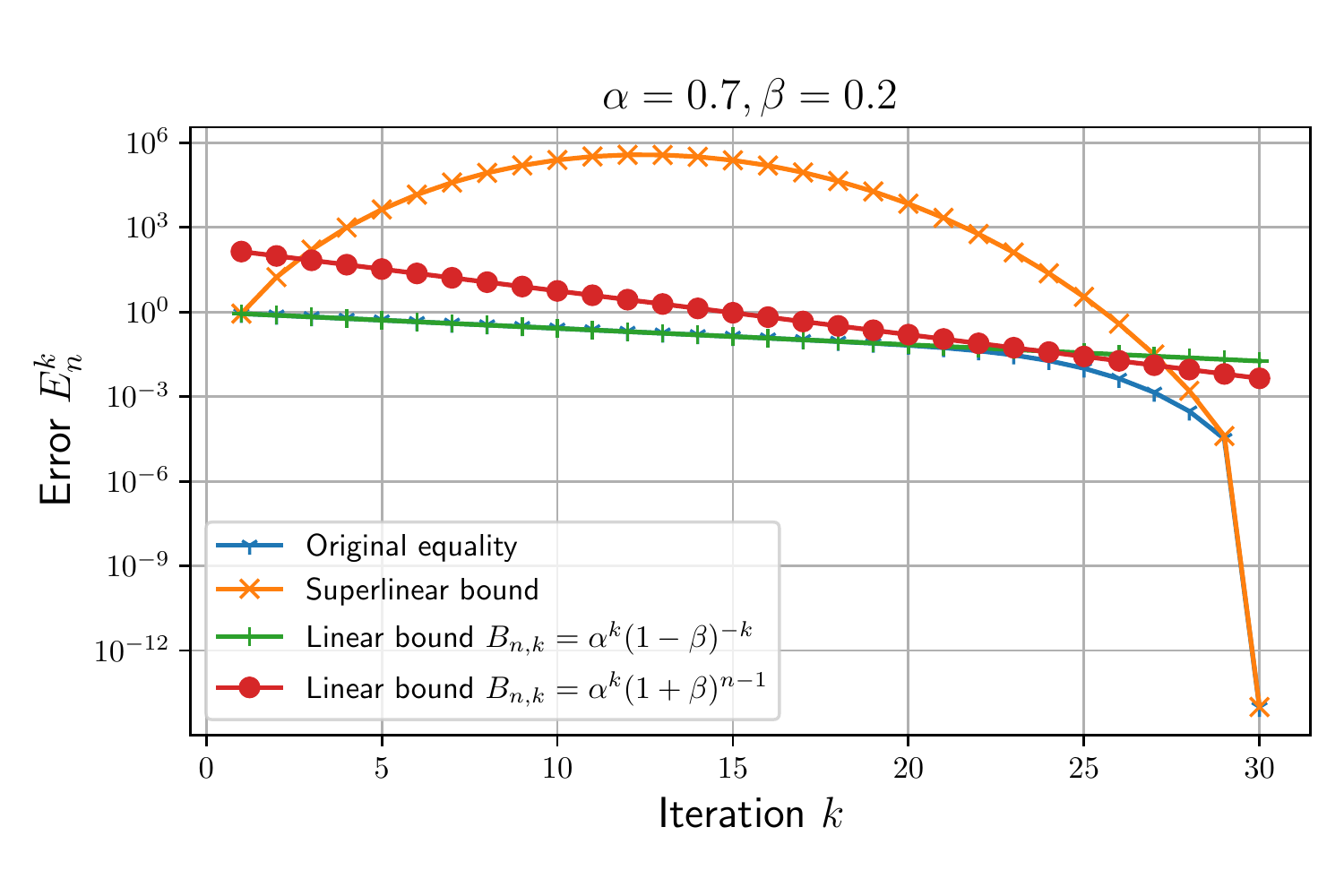}
  \end{subfigure}
  \begin{subfigure}{0.49\textwidth}
    \includegraphics[width=\textwidth]{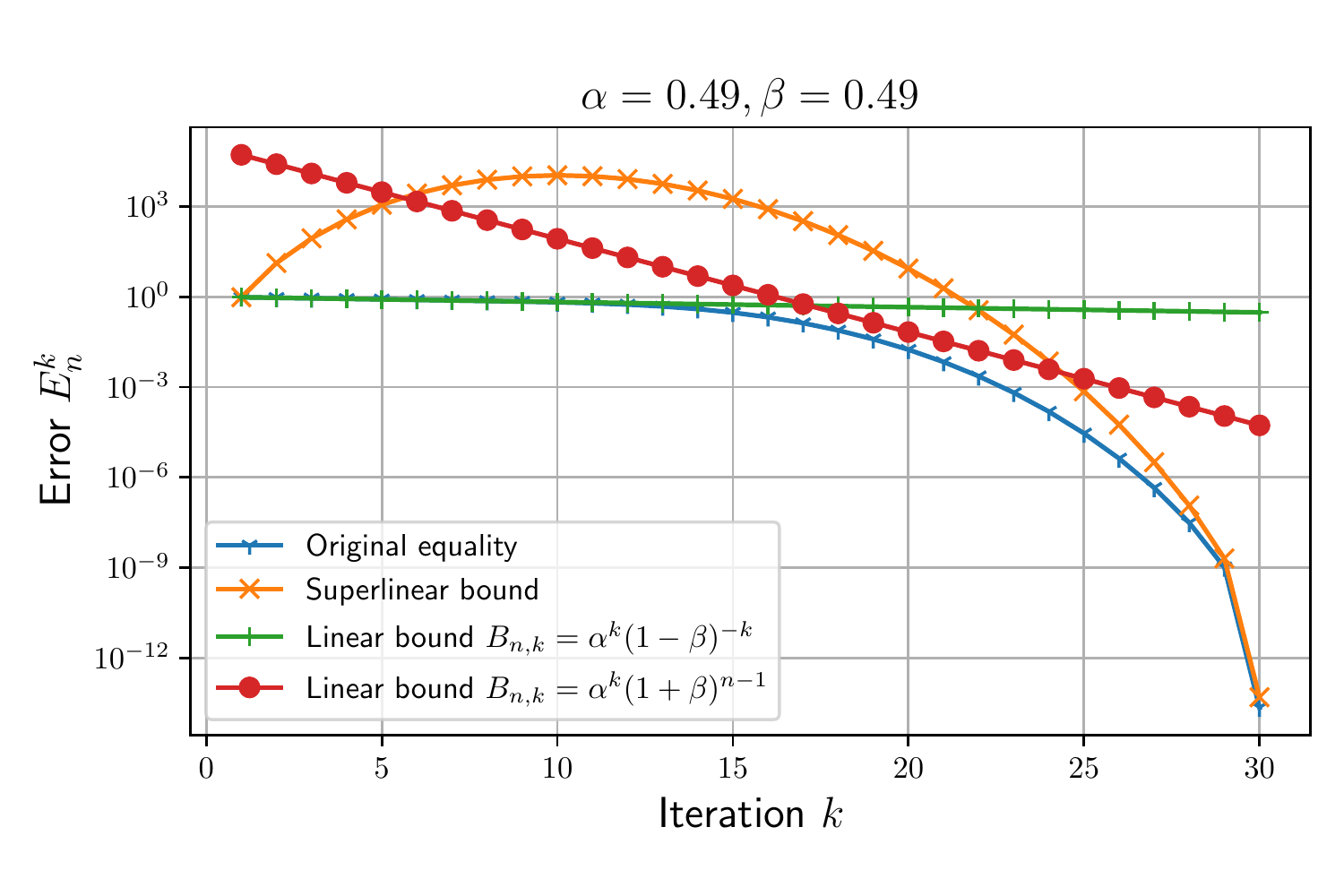}
  \end{subfigure}
  \caption{Bounds derived for several values of $\alpha$ and $\beta$. In all panels $n=30$, $\gamma=1$, and $\kappa=10^{-15}$.}
  \label{fig:bounds}
\end{figure}

In the previous section, we have proven four upper bounds for the error of low-rank Parareal. The first is directly obtained from Lemma~\ref{lem:technical_lemma}. It is the tightest bound but its expression is too unwieldy for practical use. The other three bounds can be summarized as
\begin{equation}\label{eq:summary three bounds}
  \norm{E_n^k} \leq B_{n,k} \max_{n \geq 0} \norm{E_n^0} + \frac{\kappa}{1 - \alpha - \beta},
\end{equation}
with
\begin{center}
  \begin{tabular}{ ll }\toprule
    $B_{n,k}$                                                                  & rate of~\eqref{eq:summary three bounds} in $k$ \\ \midrule
    $\alpha^k (1 - \beta)^{-k}$ \phantom{$\frac{\prod_{j=2}^k (n-j)}{(k-1)!}$} & linear                                         \\
    \addlinespace[0.5em]
    $\alpha^k (1+\beta)^{n-1}$ \phantom{$\frac{\prod_{j=2}^k (n-j)}{(k-1)!}$}  & linear                                         \\
    \addlinespace[0.5em]
    $\alpha^k (1 - \beta)^{-1} \frac{\prod_{j=2}^k (n-j)}{(k-1)!}$             & superlinear                                    \\
    \addlinespace[0.5em]
    \bottomrule
  \end{tabular}
\end{center}
Each of these practical bounds describes different phases of the convergence, and none is always better than the others. In Figure~\ref{fig:bounds}, we have plotted all four bounds for realistic values of $\alpha$ and $\beta$. We took $\kappa = 10^{-15} \approx \varepsilon_{\text{mach}}$ since it only determines the stagnation of the error and would interfere with judging the transient behavior of the convergence plot. Furthermore, the errors $e_n^0=\gamma=1$ at the start of the iteration $k=0$ were chosen arbitrarily since they have little influence on the results.

The bounds above depend on $\alpha=e^{\ell h} C_{r,q}$ and $\beta=e^{\ell h} C_q$, where $C_q$ and $C_{r,q}$ are the Lipschitz constants of $\Tq$ and $\mathcal T_{r,q}$ respectively; see~\eqref{eq:def constants abck}. While it seems difficult to give a priori results on the size $C_q$ and $C_{r,q}$, we can bound them up to first order in the theorem below. Note also that in the important case of $\ell<0$, the constants $\alpha$ and $\beta$ can be made as small as desired by taking $h$ sufficiently large.

\begin{theorem}[Lipschitz constants] \label{thm:Lipschitz_truncation}
  Let $A , \tilde{A} \in \R^{m \times n}$. Then
  \begin{align} \label{eq:Lipschitz_truncation}
    \normStd{\Tq(A) - \Tq(\tilde{A})} \leq \frac{\sigma_q}{\sigma_q - \sigma_{q+1}} \normStd{A - \tilde{A}} + O(\normStd{A - \tilde{A}}^2),
  \end{align}
  where $\sigma_q$ is the $q$th singular value of $A$. Moreover,
  \begin{align} \label{eq:Lipschitz_truncation_second}
    \normStd{\mathcal T_{r,q} (A) - \mathcal T_{r,q} (\tilde{A})} \leq \left( \frac{\sigma_q}{\sigma_q - \sigma_{q+1}} + \frac{\sigma_r}{\sigma_r - \sigma_{r+1}} \right) \normStd{A - \tilde{A}} + O(\normStd{A - \tilde{A}}^2).
  \end{align}
\end{theorem}

\begin{proof}
  For the first inequality, we refer to~\cite[Theorem 2]{breiding_sensitivity_2021} and~\cite[Theorem 24]{feppon_geometric_2018}.
  The second inequality follows from the first by the triangle inequality,
  \begin{align*}
    \normStd{\mathcal T_{r,q} (A) - \mathcal T_{r,q} (\tilde{A})}_F
     & \leq \normStd{\mathcal T_q(A) - \mathcal T_q(\tilde{A})} + \normStd{\mathcal T_r(A) - \mathcal T_r(\tilde{A})}.  \qed
  \end{align*}
\end{proof}

In many applications with low-rank matrices, the singular values of the underlying matrix are rapidly decaying.
In particular, when the singular values decay exponentially like $\sigma_k \approx e^{-ck}$ for some $c > 0$, we have
\begin{equation} \label{eq:singular gap}
  \frac{\sigma_q}{\sigma_q - \sigma_{q+1}} = \frac{1}{1- \sigma_{q+1}/\sigma_{q}} \approx \frac{1}{1 - e^{-c}}.
\end{equation}
This last quantity decreases quickly to $1$ when $c$ grows. Even for $c=1$, it is less than $1.6$. We therefore see that the constants in Theorem~\ref{thm:Lipschitz_truncation} are not too large in this case.

\begin{remark}
  In the analysis, a sufficiently large gap in the singular values is required at both the coarse rank and the fine rank. In our experiments, we observed that such a gap is indeed required at the coarse rank, but not at the fine rank. It suggests that the bound \eqref{eq:Lipschitz_truncation_second} can therefore probably be improved.
\end{remark}

\section{Numerical experiments} \label{sec:numerical_experiments}

We now show numerical experiments for our low-rank Parareal algorithm. We implemented the algorithm in Python 3.10 and all computations were performed on a MacBook Pro with a M1 processor and 16GB of RAM. The complete code is available at \href{https://github.com/BenjaminCarrel/Low-rank-Parareal}{GitHub} so that all the experiments can be reproduced. The DLRA steps are solved by the second-order projector-splitting integrator from~\cite{lubich_projector-splitting_2014}. Since the problems considered are stiff, we used sufficiently many substeps of this integrator so that the coarse and fine solvers within low-rank Parareal can be considered exact.

\subsection{Lyapunov equation} \label{sec:lyapunov}
Consider the differential Lyapunov equation,
\begin{equation} \label{eq:Lyapunov}
  \dt{X}(t) = A X(t) + X(t) A + C C^T, \quad X(0) = X_0,
\end{equation}
where $A \in \R^{m \times m}$ is a symmetric matrix, and $C \in \R^{m \times k}$ is a tall matrix for some $k \leq m$. This initial value problem admits a unique solution for $t \in [0,T]$ for any $T>0$. The most typical example of~\eqref{eq:Lyapunov} is the heat equation on a square with separable source term. Other applications can be found in~\cite{menaNumericalLowrankApproximation2018a}.

\begin{assumption} \label{ass:A is negative definite}
  The matrix $A \in \R^{m \times m}$ is symmetric and strictly negative definite.
\end{assumption}
Under Assumption~\ref{ass:A is negative definite}, the one-sided Lipschitz constant $\ell$ for~\eqref{eq:Lyapunov} is strictly negative. Indeed, the linear Lyapunov operator $\mathcal A(X) =  AX + XA$ has the symmetric matrix representation $A \otimes I + I \otimes A$ with eigenvalues $\lambda_i(A) + \lambda_j(A)$ for $1 \leq i,j \leq m$; see~\cite[Ch.~12.3]{golub_matrix_2013}. As in~\cite[Ch.~I.10]{hairer_solving_1987}, we therefore get immediately that $\ell = 2 \max_i \lambda_i(A) < 0$. Moreover, since $\mathcal A$ is invertible, we can write the closed-form solution of~\eqref{eq:Lyapunov} as
\begin{align} \label{eq:closed_form_Lyapunov}
  X(t) = e^{t \mathcal A} (X_0) + \mathcal A^{-1} (e^{t \mathcal A } (CC^T) - CC^T),
\end{align}
which can be easily verified by differentiation using properties of the matrix exponential $e^{t \mathcal A}(Z) = e^{t A} Z e^{t A}$.

The following result shows that the solution of~\eqref{eq:Lyapunov} can be well approximated by low rank. It is the analogue to a similar result for the algebraic Lyapunov equation $\mathcal{A}(X)=CC^T$. The latter result is well known, but we did not find a proof for the former in the literature.

\begin{lemma}[Low-rank approximability of Lyapunov ODE] \label{lem:low_rank_conservation}
  Let $\sigma_i(X_0)$ be the $i$th singular value of $X_0$ and likewise for $\sigma_i(CC^T)$. Under Assumption~\ref{ass:A is negative definite}, the solution $X(t)$ of~\eqref{eq:Lyapunov} has an approximation
  \[
    \text{$Y(t)$ of rank at most $r_0 + 2 r \rho$}
  \]
  for any $0 \leq r_0, r, \rho \leq m$ with error
  \[
    \norm{X(t) - Y(t)}_2 \leq e^{\ell t} \sigma_{r_0+1}(X_0) + \left( \frac{e^{t \ell} - 1}{\ell} \right) \left( 4 \exp \left( \frac{- \pi^2 \rho}{\log(4 \kappa_A)} \right) \| CC^T \|_2 + \sigma_{r+1}(CC^T) \right),
  \]
  where $\kappa_A = \norm{A}_2 \| A^{-1}\|_2$ and $\ell = 2 \max_i \lambda_i(A)$.

\end{lemma}

\begin{proof}
  The aim is to approximate the following two terms that make up the closed-form solution $X(t)$ in~\eqref{eq:closed_form_Lyapunov}:
  \begin{align*}
    X_1(t) = e^{t \mathcal A} (X_0), \quad X_2(t) = \mathcal A^{-1} (e^{t \mathcal A } (CC^T) - CC^T).
  \end{align*}
  The first term $X_1(t)$ can be treated directly. By the truncated SVD, the initial value satisfies
  $$X_0 = Y_0 + E_0 \text{ where } \rank(Y_0) = r_0 \text{ and } \|E_0\|_2 = \sigma_{r_0+1}(X_0).$$
  By Assumption \ref{ass:A is negative definite}, the operator $\mathcal A$ is full rank. We therefore obtain
  \begin{align} \label{eq:first term Lyapunov}
    X_1(t) = e^{t \mathcal A} (X_0) =  e^{At} Y_0 e^{At} + e^{At} E_0 e^{At} = Y_1(t) + E_1(t),
  \end{align}
  where $\rank (Y_1(t)) = \rank (Y_0) = r_0$ and $\| E_1(t)\|_2 \leq e^{\ell t} \sigma_{r_0+1}(X_0)$ since $\ell = 2 \max_i \lambda_i(A)$.

  Next, we focus on the second term $X_2(t)$. Like above, the source term can be decomposed as
  $$CC^T = D + F \text{ where } \rank(D) = r \text{ and } \| F \|_2 = \sigma_{r+1}(CC^T).$$
  By linearity of the Lyapunov operator, we therefore obtain
  \begin{align} \label{eq:second term Lyapunov}
    X_2(t)
    = \mathcal A^{-1} (e^{t \mathcal A} D - D) + \mathcal{A}^{-1} (e^{t \mathcal A} F - F).
  \end{align}
  Denote $M=e^{t \mathcal A} D - D$. By definition of the Lyapunov operator $\mathcal{A}$, we have
  \begin{align*}
    S = \mathcal A^{-1} (M) \iff AS + SA = M.
  \end{align*}
  As studied in~\cite{penzl_eigenvalue_2000} and then improved in~\cite{beckermann_singular_2017}, the singular values of the solution $S$ are bounded as
  \begin{equation}\label{eq:bound_decay_Lyap}
    \frac{ \sigma_{\rank(M) \rho + 1}(S) }{\sigma_1(S)} \leq 4 \exp \left( \frac{- \pi^2 \rho}{\log(4 \kappa_A)} \right),
  \end{equation}
  where $\kappa_A = \norm{A}_2 \| A^{-1}\|_2$ and $ 0 \le \rho \le m$. Since $\rank(M) \leq 2 \rank(D) = 2r$ by assumption on $D$, the bound~\eqref{eq:bound_decay_Lyap} then implies that
  \[
    S = Y_2(t) + \delta S(t),
  \]
  where $\rank(Y_2(t)) \leq 2 r \rho$ and
  \begin{align*}
     & \norm{\delta S(t)}_2 \leq 4 \exp \left( \frac{- \pi^2 \rho}{\log(4 \kappa_A)} \right) \norm{\mathcal A^{-1}(e^{t \mathcal A} D - D)}_2 \leq 4 \exp \left( \frac{- \pi^2 \rho}{\log(4 \kappa_A)} \right) \frac{e^{t \ell} - 1}{\ell} \norm{D}_2,
  \end{align*}
  where the last inequality holds by properties of the logarithmic norm $\mu$ of $\mathcal{A}$ which equals $\ell$; see \cite[Proposition 2.2]{soderlind_logarithmic_2006}.
  Moreover, we can bound the last term in \eqref{eq:second term Lyapunov} as
  \begin{align*}
    \| E_2(t) \|_2 = \| \mathcal{A}^{-1} (e^{t \mathcal A} F - F) \|_2 \leq  \frac{e^{t \ell} - 1}{\ell} \norm{F}_2.
  \end{align*}
  Putting \eqref{eq:first term Lyapunov} and \eqref{eq:second term Lyapunov} together, we obtained
  \begin{align*}
    X(t) = Y(t) + E(t), \quad Y(t) = Y_1(t) + Y_2(t), \quad E(t) = E_1(t) + \delta S(t) + E_2(t),
  \end{align*}
  which proves the statement of the lemma.
\end{proof}

The lemma shows that if $X_0$ and $CC^T$ have good low-rank approximations, then the solution $X(t)$ of the differential Lyapunov equation has comparable low-rank approximations as well on $[0,T]$. Since $\ell < 0$, we can even take $T \to \infty$ and recover essentially the low-rank approximability of the Lyapunov equation $X(\infty) = \mathcal{A}^{-1}(CC^T)$. This is clearly visible when $X_0$ and $CC^T$ are exactly of low rank, which we state as a simple corollary for convenience.

\begin{corollary}\label{cor:low-rank-approx-diff-Lyap}
  Under Assumption \ref{ass:A is negative definite} and assuming that $\rank(X_0) = r_0, \, \rank(CC^T) = r,$ the solution $X(t)$ of \eqref{eq:closed_form_Lyapunov} has an approximation
  \[
    \text{$Y(t)$ of rank at most $r_0 + 2 r \rho$}
  \]
  for any $0 \le \rho \leq m$ with error
  \begin{align*}
    \norm{X(t) - Y(t)}_2 \leq 4\frac{e^{t \ell} - 1}{\ell} \exp \left( \frac{- \pi^2 \rho}{\log(4 \kappa_A)} \right) \norm{CC^T}_2.
  \end{align*}

\end{corollary}

The corollary clearly shows that the approximation error decreases exponentially when the approximation rank increases linearly via $\rho$. Furthermore, we see that the condition number of the matrix $A$ has only a mild influence due to $\log(\kappa_A)$.

\begin{remark}
  Corollary~\ref{cor:low-rank-approx-diff-Lyap} can be compared to a similar result in~\cite{koskela_analysis_2020}. In that work, the authors solve~\eqref{eq:Lyapunov} with exact low-rank $X_0=ZZ^T$ and $CC^T$ using a Krylov subspace method. More specifically, with $U_k$ an orthonormal matrix that spans the block Krylov space $K_k(A, [C\ Z])$, the projected Lyapunov equation
  $$
    \dt{Y_k} = (U_k^T A U_k)\, Y_k + Y_k \, (U_k^T A U_k) + U_k^T C C^T U_k, \quad Y_k(0) = U_k^T X_0 U_k,
  $$
  is used to define the approximation $X_k(t) = U_k Y_k(t) U_k^T$. The approximation error of $X_k(t)$ is studied in~\cite[Theorem 4.2]{koskela_analysis_2020}. Since $\rank(X_k(t)) \leq k(\rank(Z) + \rank(C))$, we therefore also get a result on the low-rank approximability of~\eqref{eq:Lyapunov}. This bound is, however, worse than ours since it does not give zero error for $t=0$ and $k=1$, for example. On the other hand, it is a bound for a discrete method whereas our Lemma~\ref{lem:low_rank_conservation} and Corollary~\ref{cor:low-rank-approx-diff-Lyap} are statements about the exact solution.
\end{remark}

We now apply the low-rank Parareal algorithm to the differential Lyapunov equation~\eqref{eq:Lyapunov}. Let $A = \Delta_{dx}$ be the $n \times n$ discrete Laplacian with zero Dirichlet boundary conditions obtained by standard centered differences on $[-1, 1]$. The Lyapunov equation is therefore a model for the 2D heat equation on $\Omega = [-1,1]^2$. In the experiments, we used $n=100$ spatial points and the time interval $[0,T] = [0,2]$. The matrix $C$ for the source is generated randomly with singular values $\sigma_i = 10^{-5(i-1)}$ where $i=1,2,\ldots$ so that its numerical rank is $4$. In order to have a realistic initial value, $X_0$ is obtained as the exact solution at time $t=0.01$ of the same ODE but with a random initial value $\tilde{X_0}$ with singular values $\sigma_i = 10^{-(i-1)}$.

\begin{figure}
  \centering
  \begin{subfigure}{0.32\textwidth}
    \includegraphics[width=\textwidth]{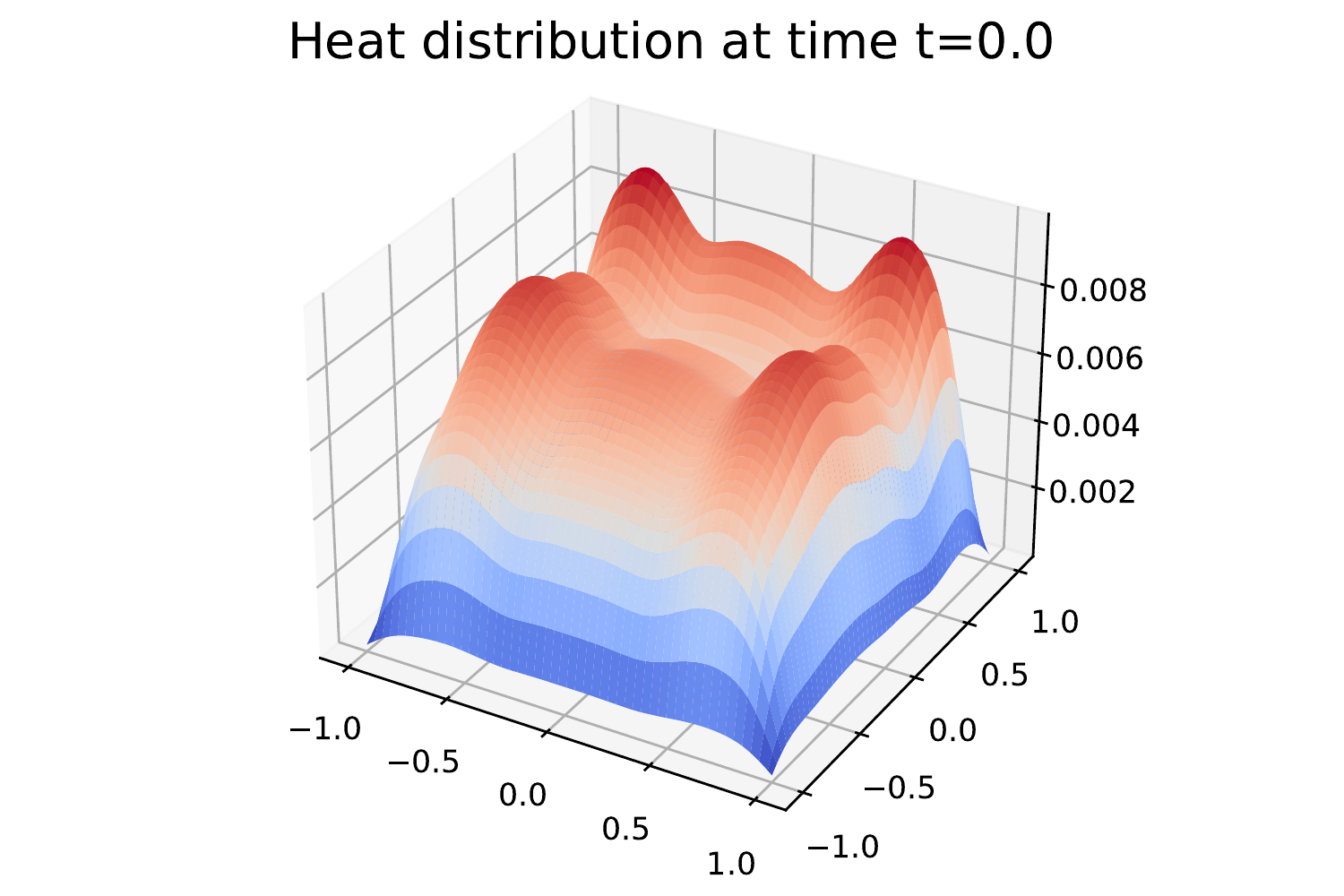}
  \end{subfigure}
  \begin{subfigure}{0.32\textwidth}
    \includegraphics[width=\textwidth]{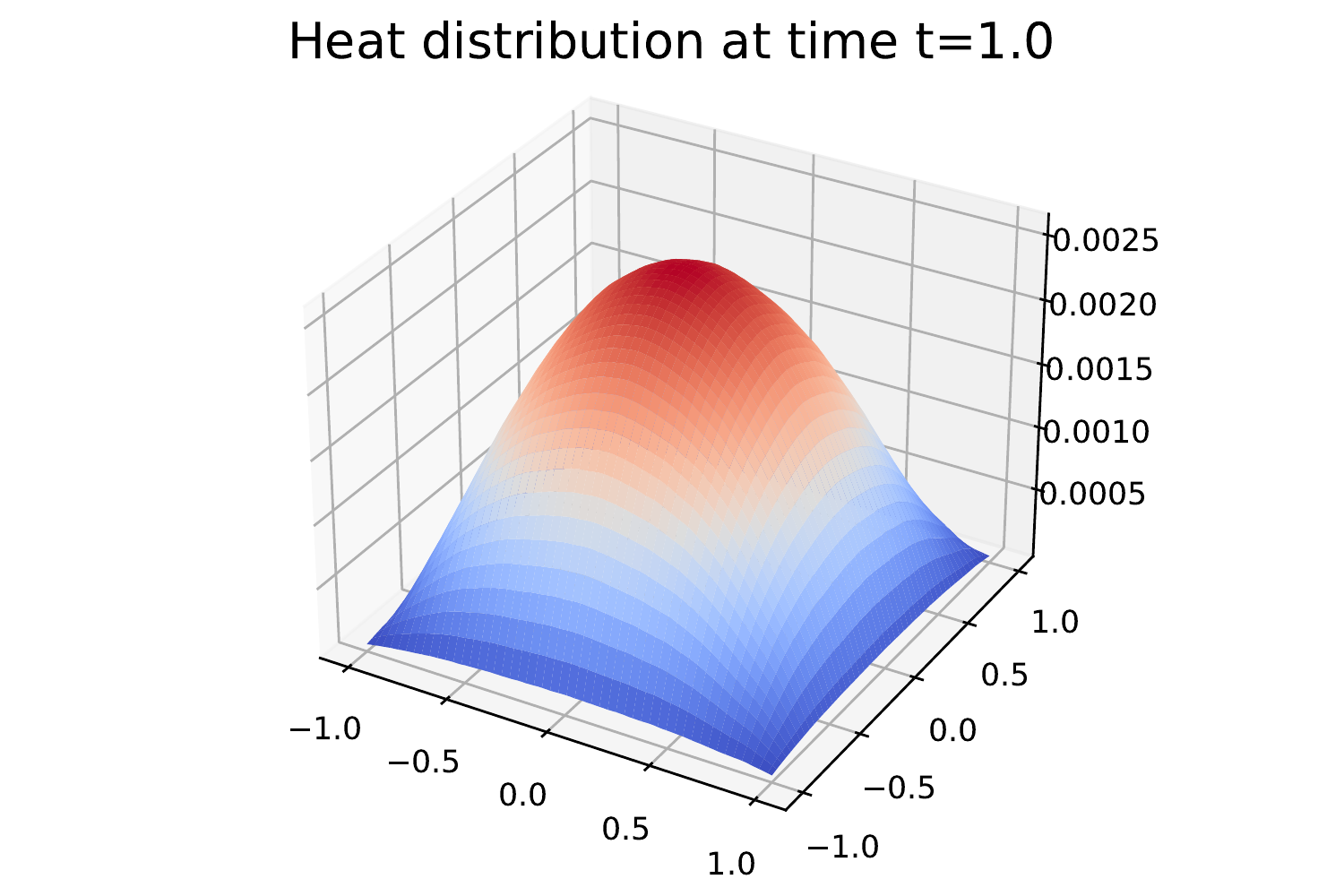}
  \end{subfigure}
  \begin{subfigure}{0.32\textwidth}
    \includegraphics[width=\textwidth]{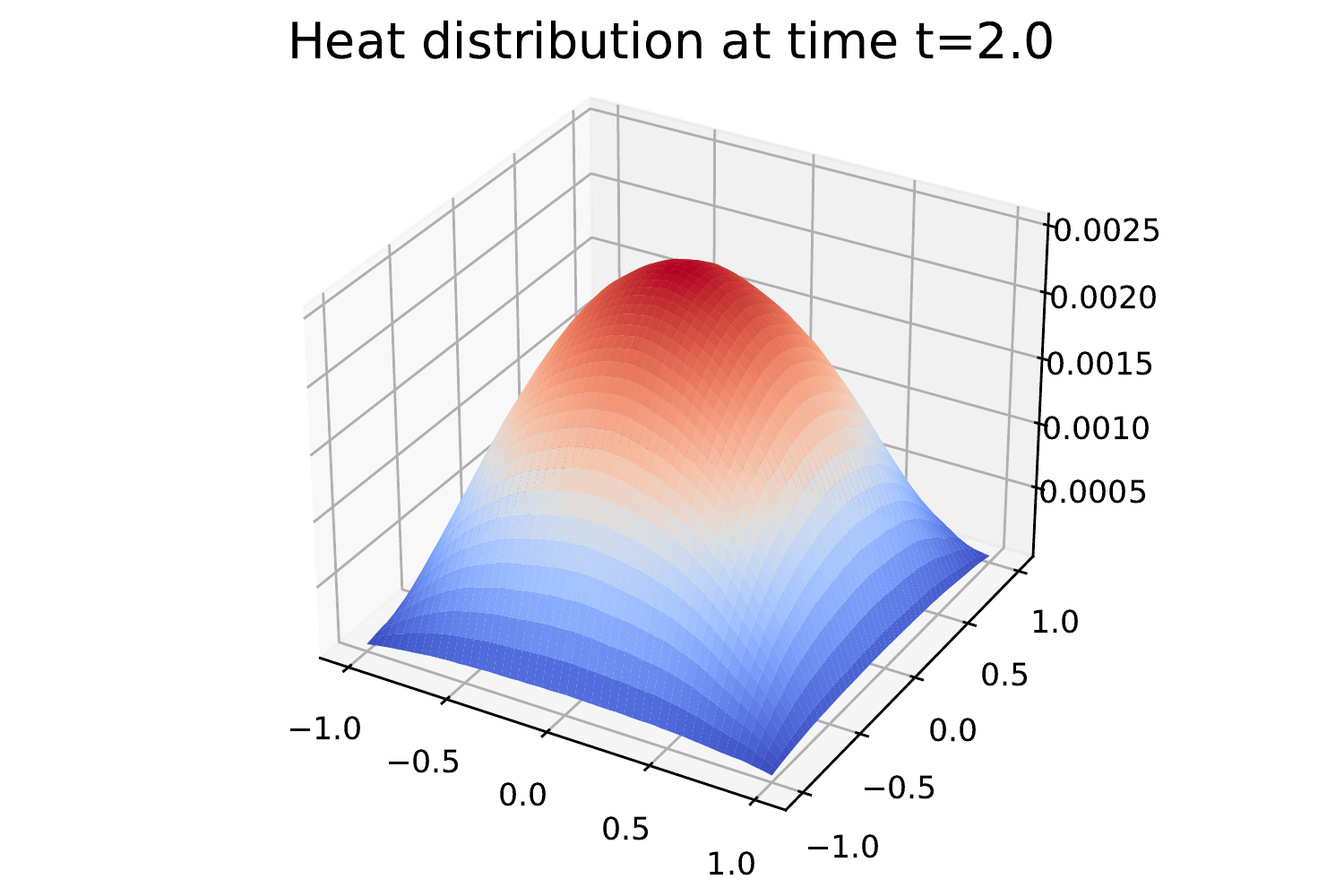}
  \end{subfigure}
  \begin{subfigure}{0.32\textwidth}
    \includegraphics[width=\textwidth]{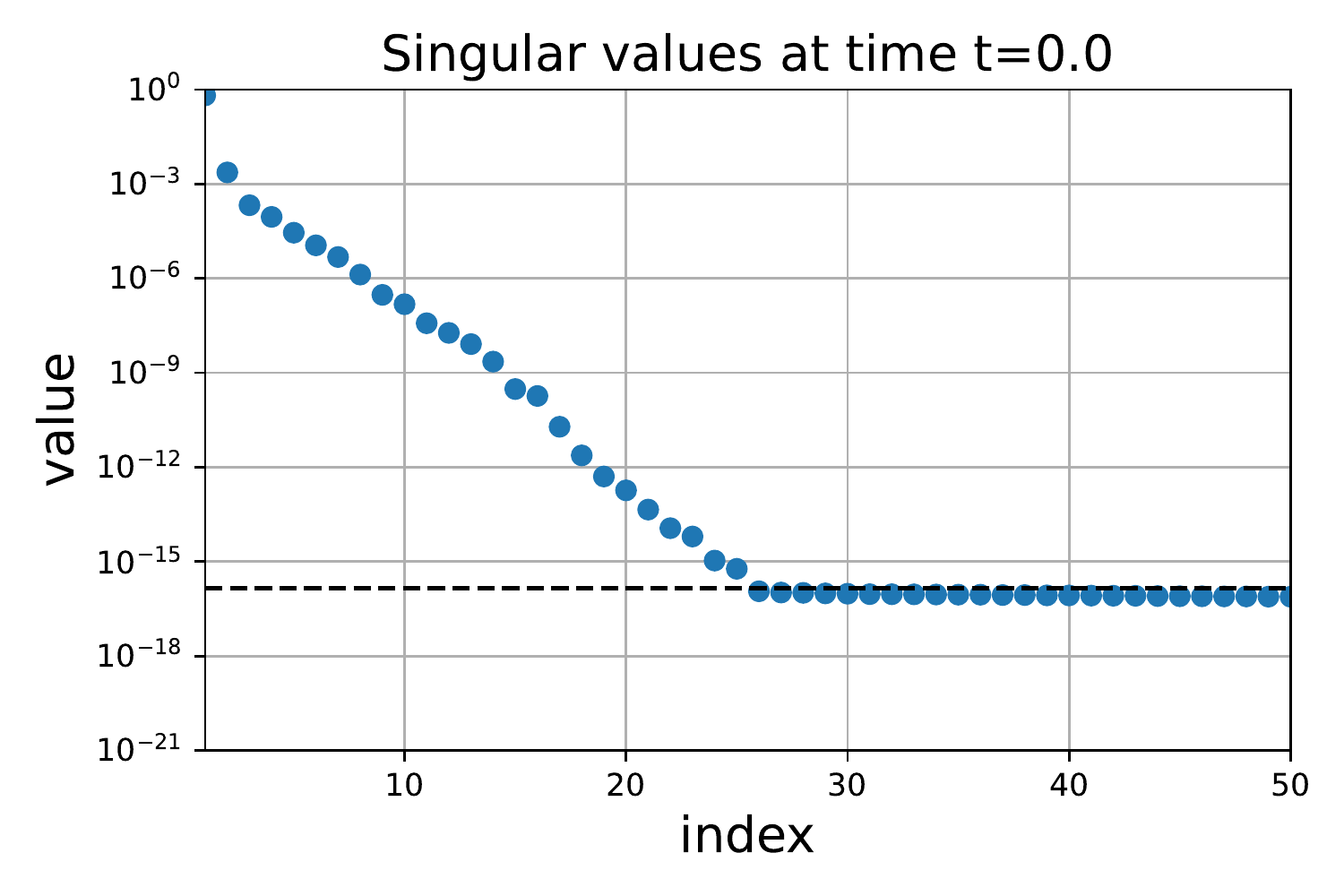}
  \end{subfigure}
  \begin{subfigure}{0.32\textwidth}
    \includegraphics[width=\textwidth]{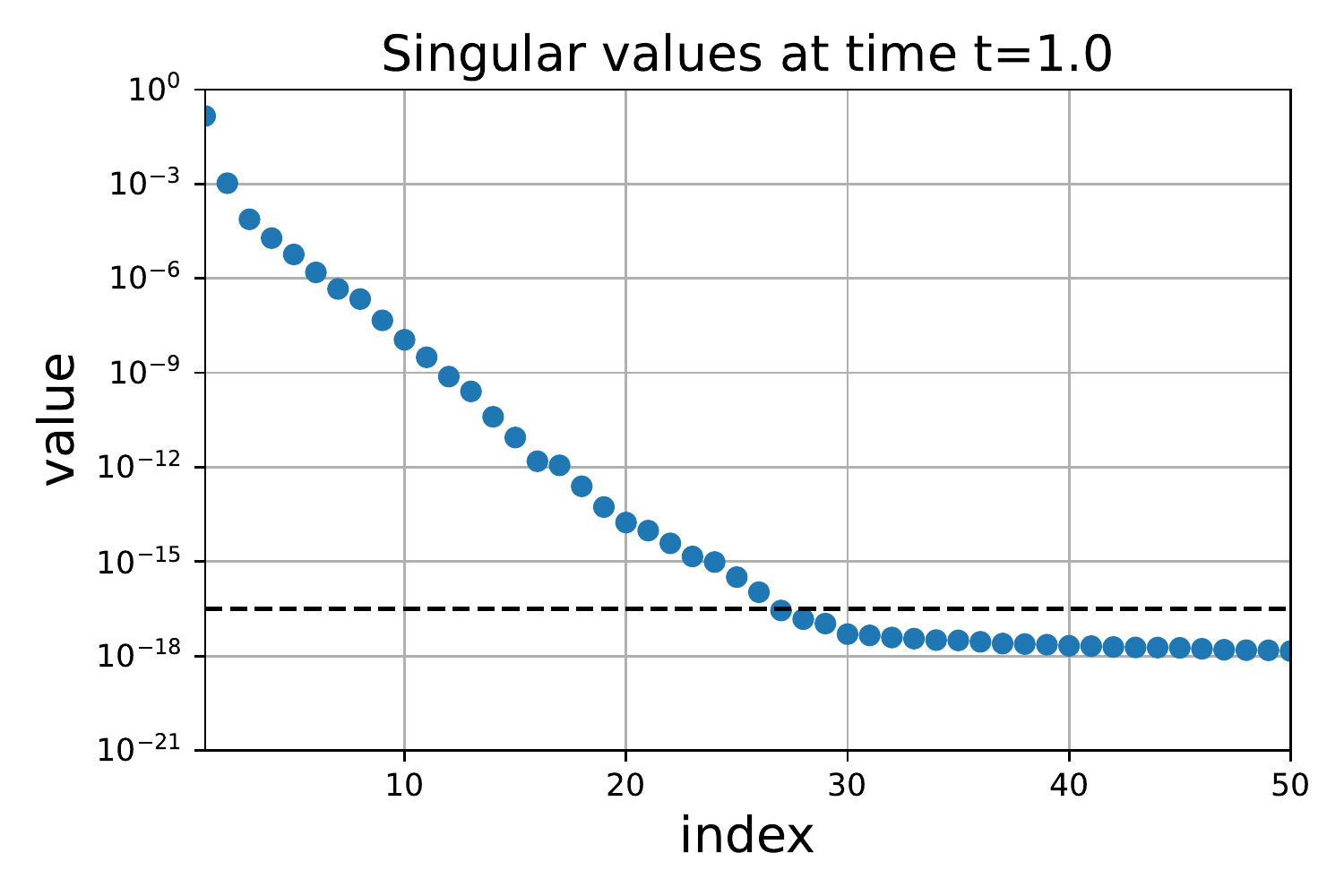}
  \end{subfigure}
  \begin{subfigure}{0.32\textwidth}
    \includegraphics[width=\textwidth]{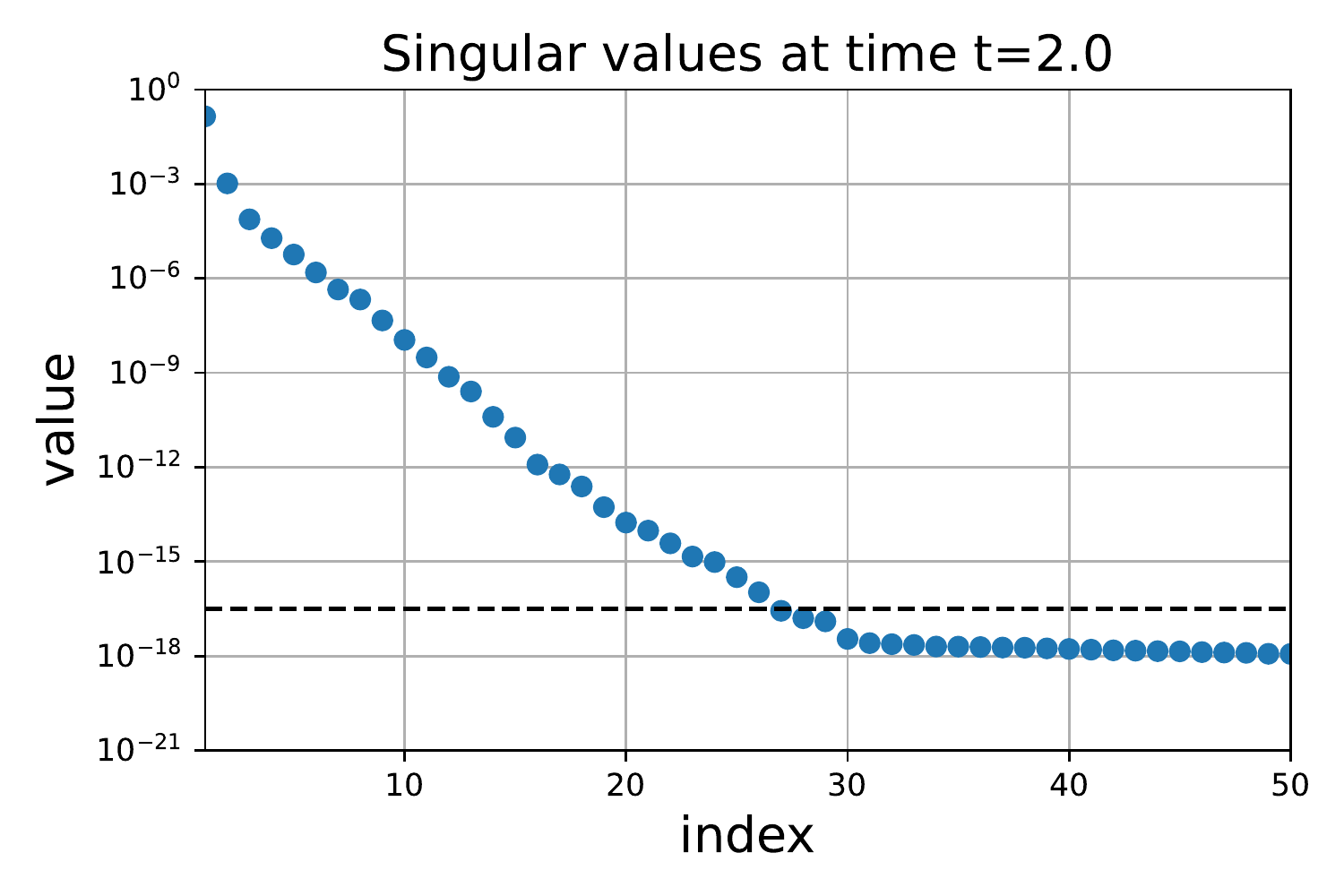}
  \end{subfigure}
  \caption{Solution over time of the Lyapunov ODE~\eqref{eq:Lyapunov} for the heat equation. Note the change of scale between $t=0.0$ and $t=1.0$.}
  \label{fig:lyapunov_solution}
\end{figure}

Figure~\ref{fig:lyapunov_solution} is a 3D plot of the solution over time on $\Omega$ with its corresponding singular values. As we can see, the solution becomes almost stationary at $t=1.0$. In addition, it stays low-rank over time in agreement to Lemma~\ref{lem:low_rank_conservation}. Moreover, the singular values suggest to take the fine rank $r=16$ for an error of the fine solver of order $10^{-12}$.

\begin{figure}
  \begin{subfigure}{0.49\textwidth}
    \includegraphics[width=\textwidth]{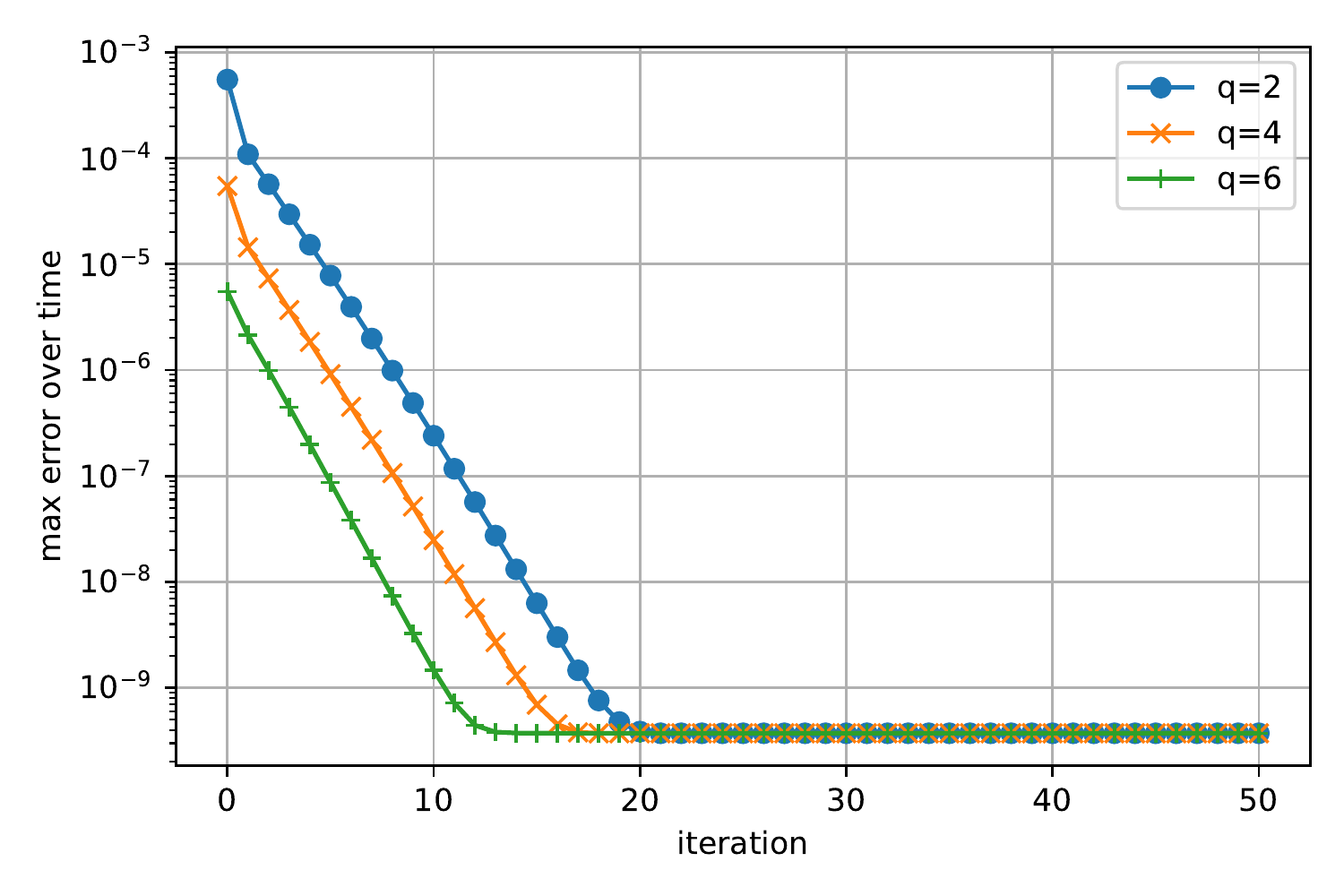}
    \caption{Several coarse ranks $q$ with fine rank $r=16$.}
    \label{fig:lyapunov_several_coarse_ranks}
  \end{subfigure}
  \begin{subfigure}{0.49\textwidth}
    \includegraphics[width=\textwidth]{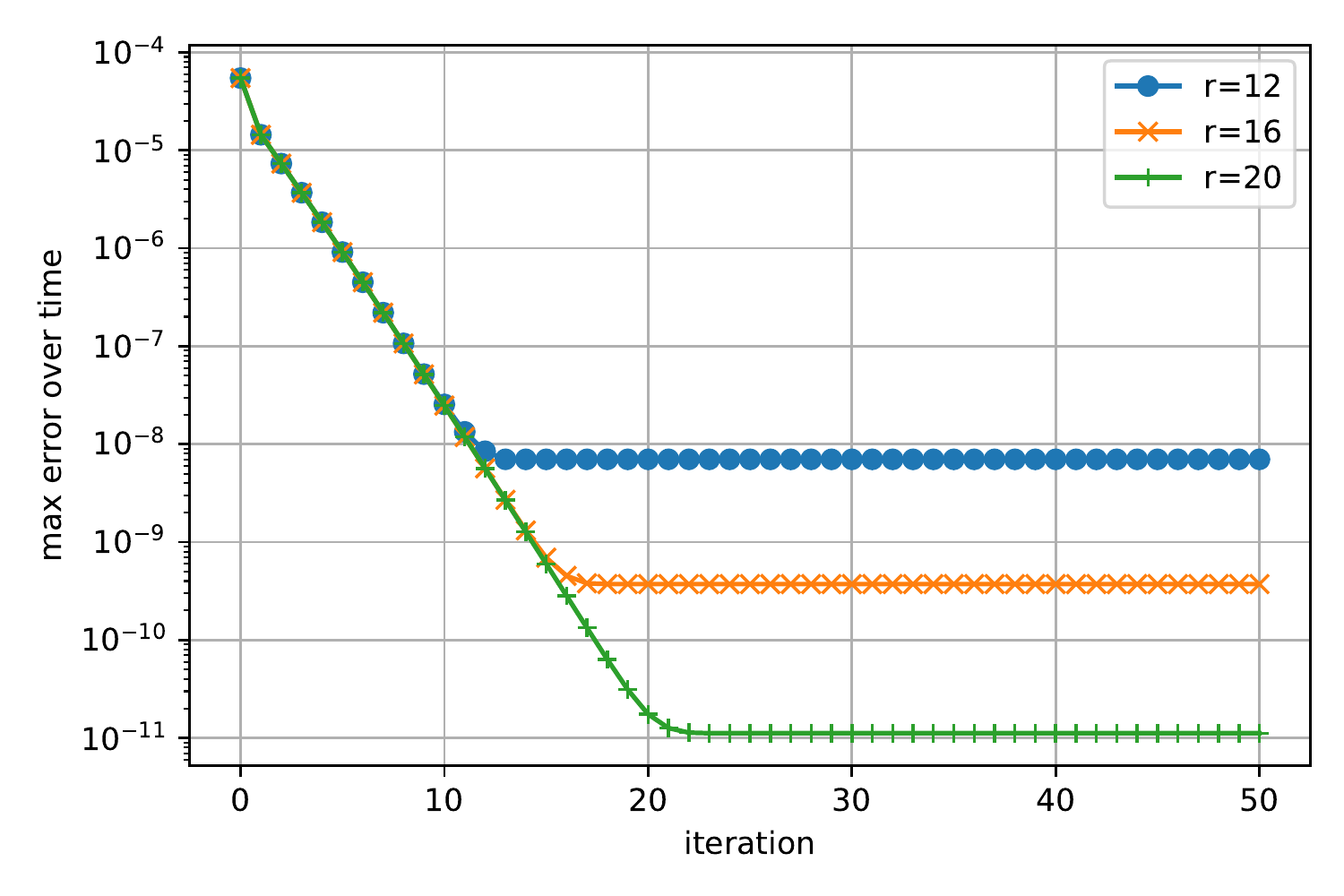}
    \caption{Several fine ranks $r$ with coarse rank $q=4$.}
    \label{fig:lyapunov_several_fine_ranks}
  \end{subfigure}
  \caption{Convergence of the error of low-rank Parareal for the Lyapunov ODE~\eqref{eq:Lyapunov} with $n=100$ and $T=2.0$. Influence of the coarse and fine ranks.}
  \label{fig:lyapunov_several_ranks}
\end{figure}

The convergence of the error of the low-rank Parareal algorithm is shown in Figure~\ref{fig:lyapunov_several_ranks}. The algorithm converges linearly from the coarse rank solution to the fine rank solution. Figure~\ref{fig:lyapunov_several_coarse_ranks} suggests that the coarse rank does not influence the convergence rate and it only reduces the initial error. This is consistent with our analysis. Indeed, since the singular values are exponentially decaying, the singular gap is approximately constant; see~\eqref{eq:singular gap}. Hence, the constants $\alpha$ and $\beta$ from \eqref{eq:def constants abck} that determine the convergence rate do not depend on the coarse rank $q$; as is shown up to first order in Theorem~\ref{thm:Lipschitz_truncation}. Figure~\ref{fig:lyapunov_several_fine_ranks} shows that, similarly, the convergence rate does not depend on the fine rank either, although it limits the final error.

\begin{figure}
  \begin{subfigure}{0.49\textwidth}
    \includegraphics[width=\textwidth]{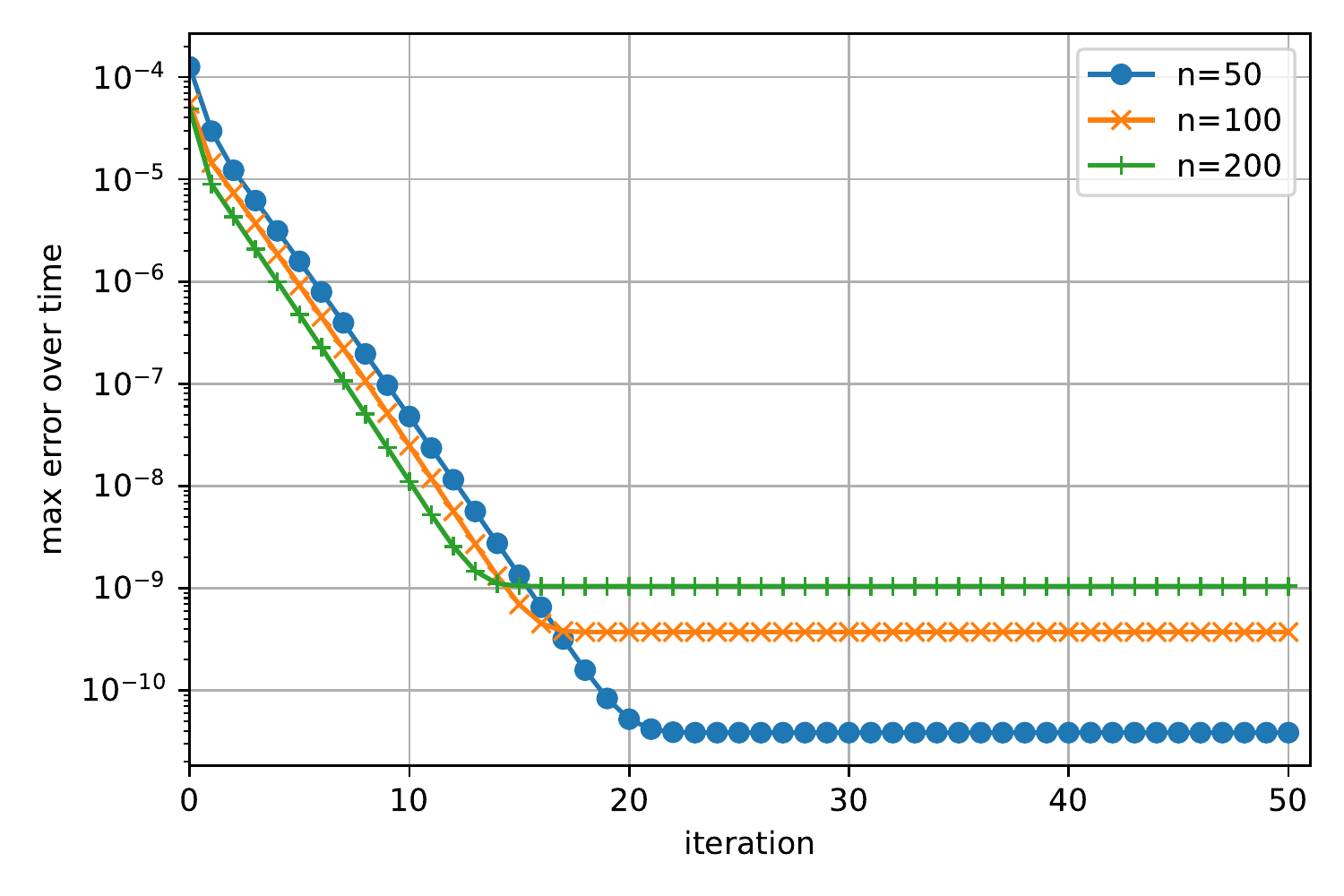}
    \caption{Several sizes $n$ with final time $T=2.0$.}
    \label{fig:lyapunov_several_sizes}
  \end{subfigure}
  \begin{subfigure}{0.49\textwidth}
    \includegraphics[width=\textwidth]{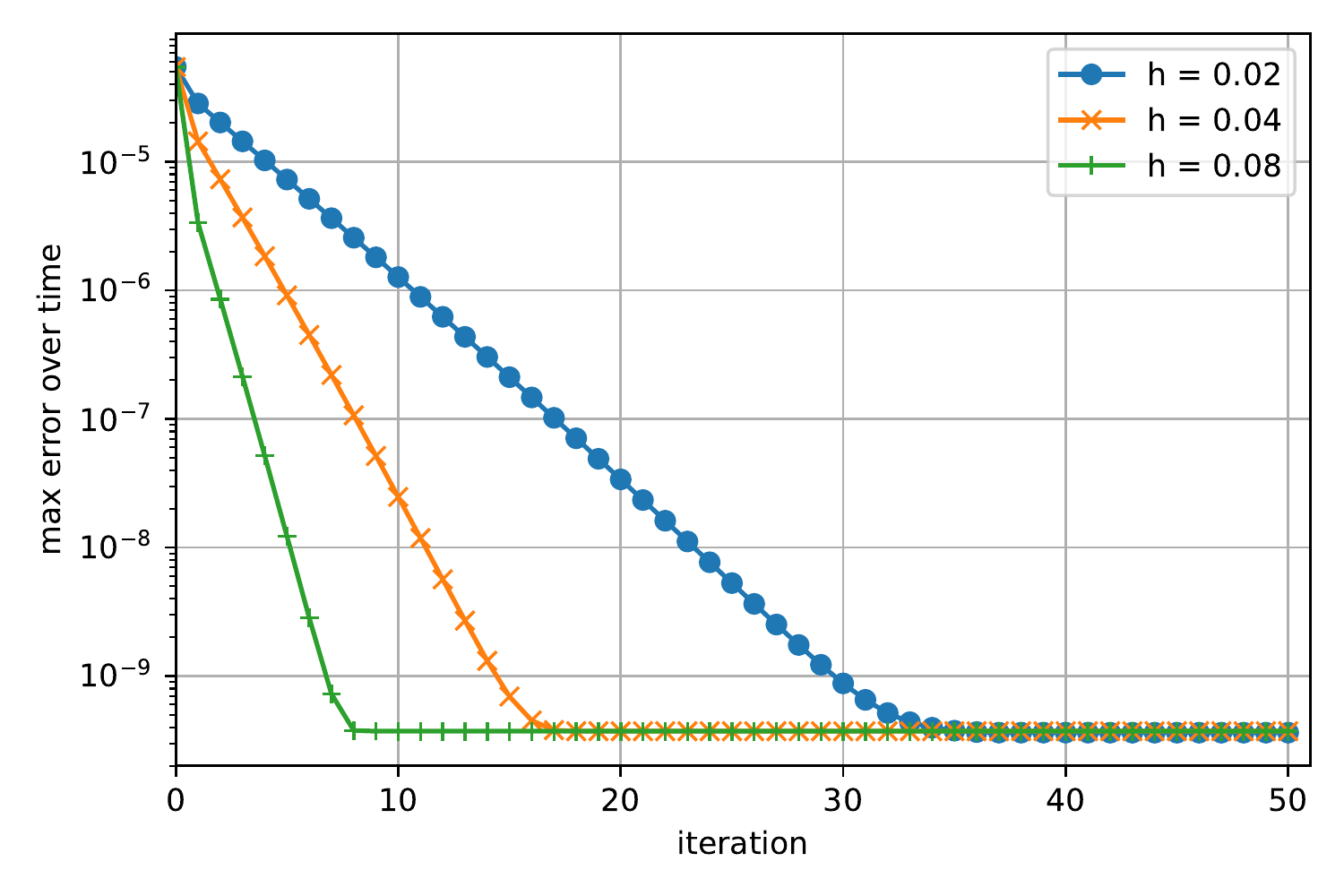}
    \caption{Several stepsizes $h$ with size $n=100$.}
    \label{fig:lyapunov_several_final_times}
  \end{subfigure}
  \caption{Convergence of the error of low-rank Parareal for the Lyapunov ODE~\eqref{eq:Lyapunov} with coarse rank $q=4$ and fine rank $r=16$. Influence of size and final time.}
\end{figure}

In Figure~\ref{fig:lyapunov_several_sizes}, we investigate the convergence for several sizes $n$. Even though the problem is stiff, the convergence does not seem influenced by the size of the problem.
Figure~\ref{fig:lyapunov_several_final_times} shows the error of the algorithm applied to the problem with several step sizes. According to our analysis, the convergence is faster when the stepsize $h$ is large; see \eqref{eq:def constants abck}.

\subsection{Parametric cookie problem} \label{sec:cookie}

We now solve a simplified version of the parametric cookie problem from~\cite{kressner_low-rank_2011}. Consider the ODE
\begin{equation} \label{eq:parametric_ode}
  \dt{Y} = -A_0 Y - A_1 Y C_1 + \mathbf{b} \mathbf{1}^T, \quad Y(0) = Y_0,
\end{equation}
where the sparse matrices $A_0, A_1 \in \R^{1580 \times 1580}$, $\mathbf{b} \in \R^{1580}$, and $C_1 = \diag(c_1^1, c_1^2, \ldots, c_1^p)$ are given in~\cite{kressner_low-rank_2011}. The aim of this problem is to solve a heat problem simultaneously with several heat coefficients, denoted by $c_1^1, \ldots, c_1^p$.

In our experiments, we used $p=101$ parameters with $c_1^1 = 0, c_1^2 = 1, \ldots, c_1^{101} = 100$. The initial value $X_0$ is obtained after computing the exact solution of~\eqref{eq:parametric_ode}  at time $t=0.01$ with the zero matrix as initial value. The time interval is $[0,T] = [0, 0.1]$.

\begin{figure}
  \centering
  \begin{subfigure}{0.32\textwidth}
    \includegraphics[width=\textwidth]{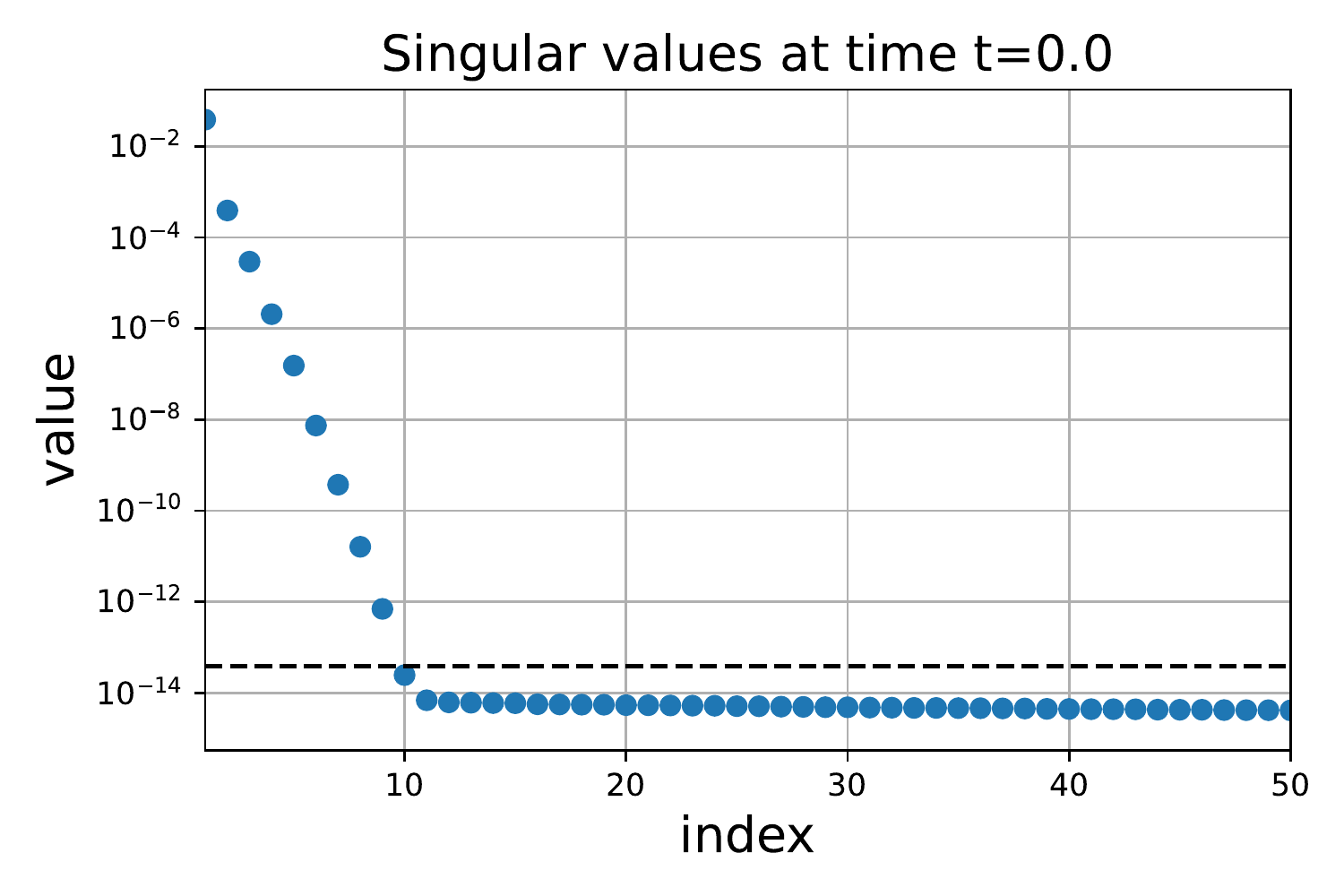}
  \end{subfigure}
  \begin{subfigure}{0.32\textwidth}
    \includegraphics[width=\textwidth]{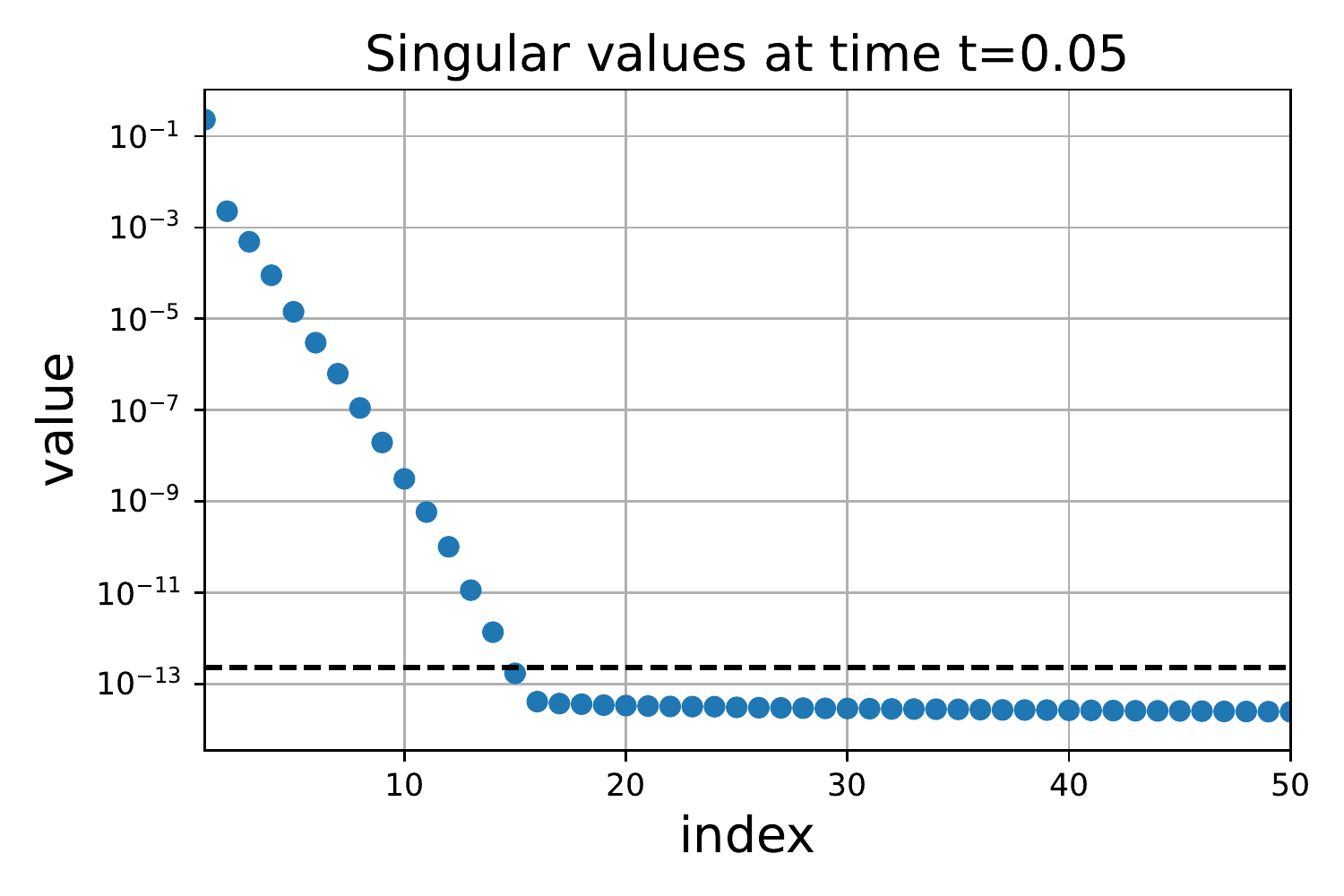}
  \end{subfigure}
  \begin{subfigure}{0.32\textwidth}
    \includegraphics[width=\textwidth]{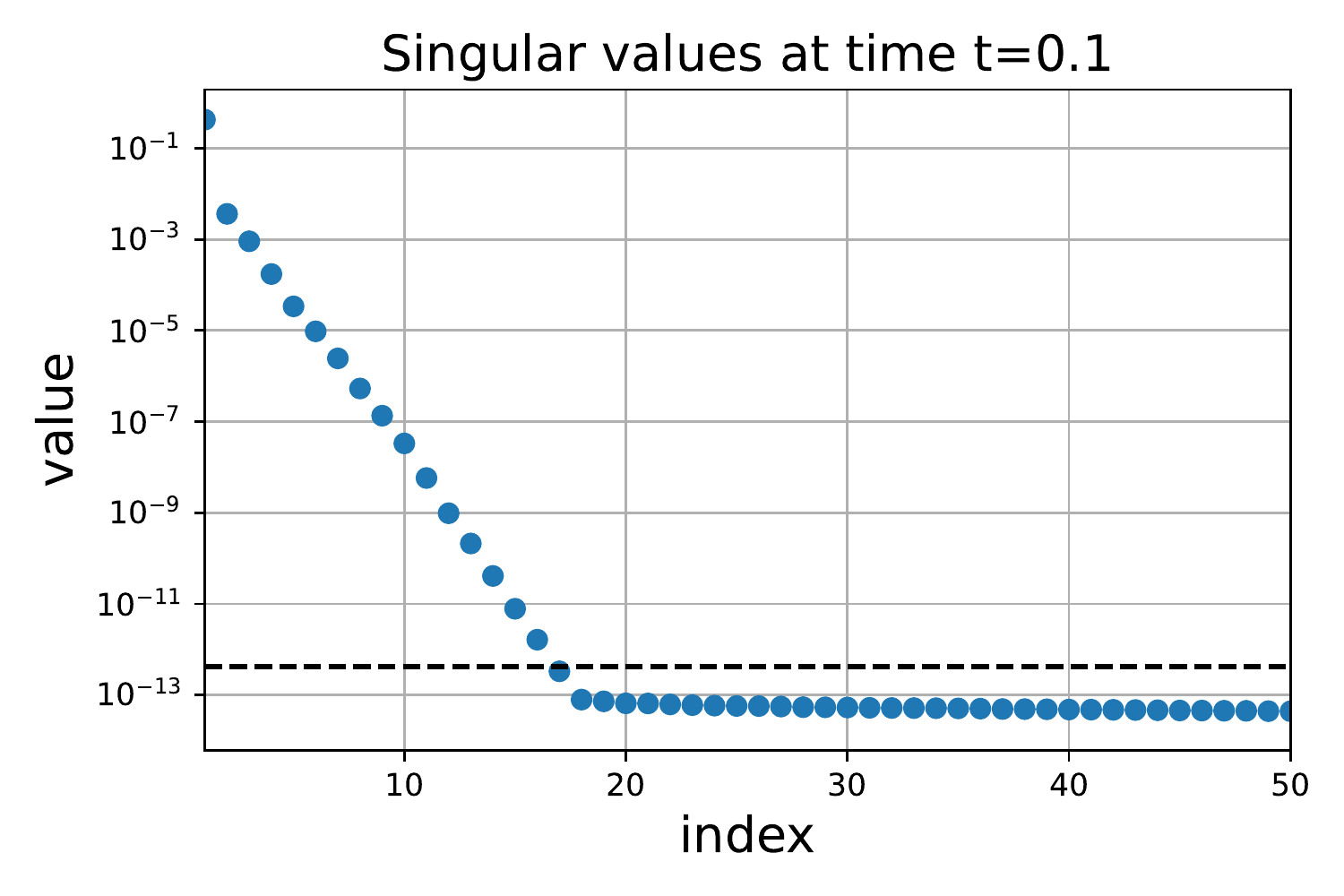}
  \end{subfigure}
  \caption{Singular values of the solution over time of the parametric cookie problem~\eqref{eq:parametric_ode}.}
  \label{fig:cookie_singular_values}
\end{figure}

The singular values of the reference solution are shown in Figure~\ref{fig:cookie_singular_values}. The stationary solution has good low-rank approximations, as was proved in~\cite[Thm.~2.4]{kressner_low-rank_2011}. The singular value decay suggests that a fine rank $r=16$ leads to full numerical accuracy.

\begin{figure}
  \centering
  \begin{subfigure}{0.49\textwidth}
    \includegraphics[width=\textwidth]{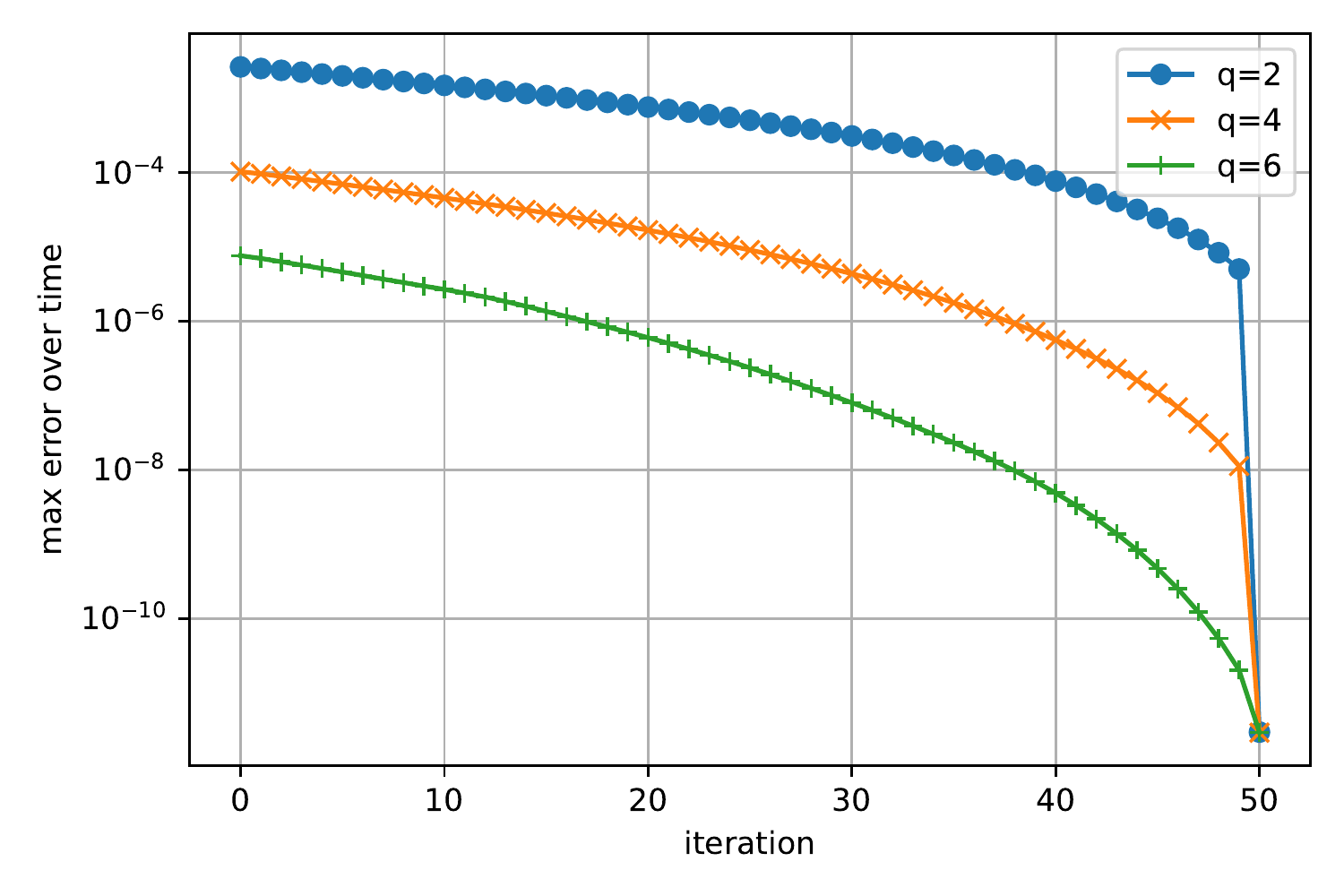}
    \caption{Several coarse ranks with fine rank $r=16$.}
  \end{subfigure}
  \begin{subfigure}{0.49\textwidth}
    \includegraphics[width=\textwidth]{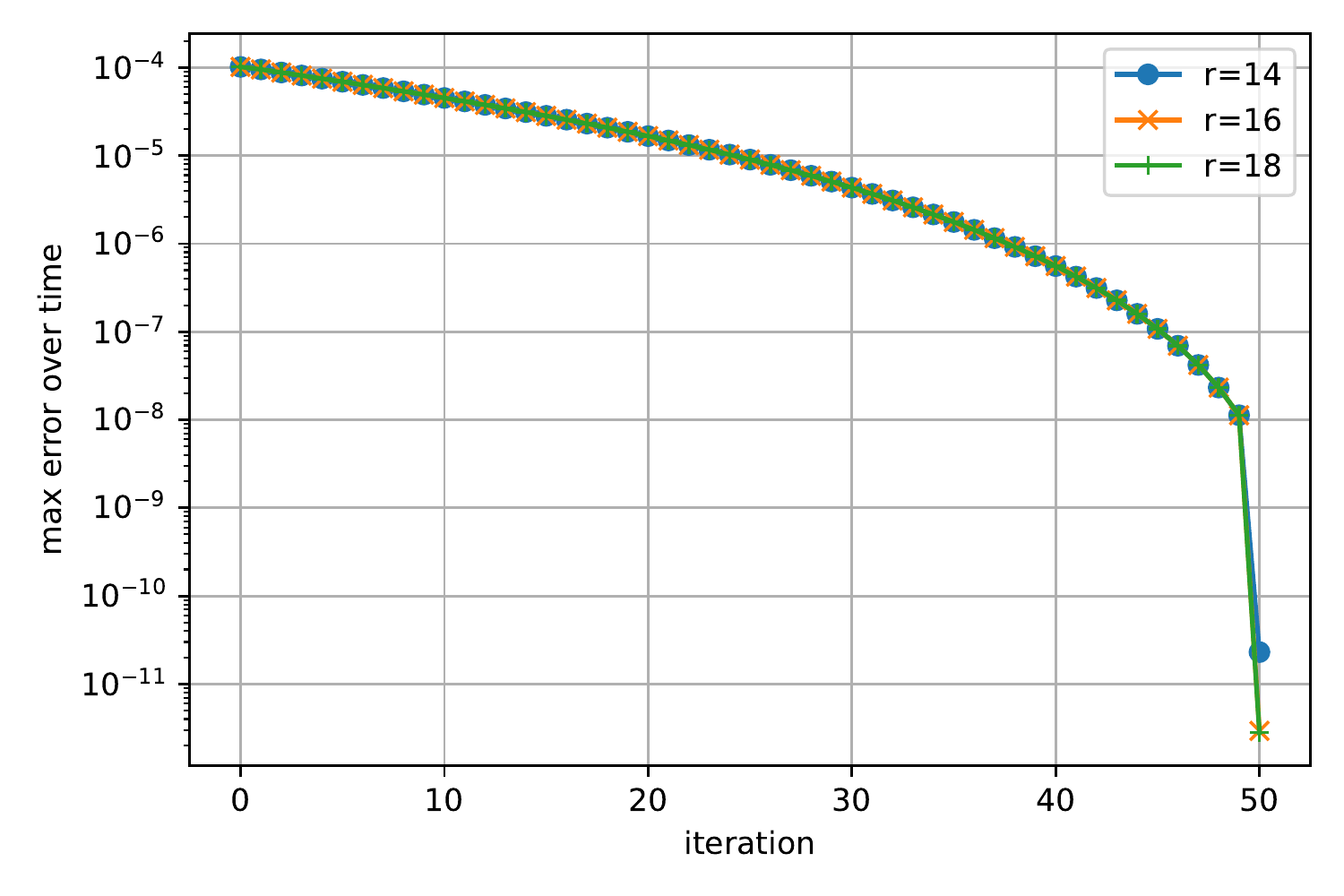}
    \caption{Several fine ranks with coarse rank $q=4$.}
  \end{subfigure}
  \caption{Convergence of the error of low-rank Parareal for the parametric cookie problem~\eqref{eq:parametric_ode}. Influence of the coarse and fine ranks.}
  \label{fig:cookie_several_ranks}
\end{figure}

In Figure~\ref{fig:cookie_several_ranks}, we applied the low-rank Parareal algorithm with several coarse ranks $q$ and fine ranks $r$. Like for the Lyapunov equation, it seems that the convergence rate does not depend on the coarse rank $q$. In agreement to our analysis (see Figure~\ref{fig:bounds}), the convergence is linear in the first iterations and superlinear in the last iterations. In addition, the convergence is not influenced by the fine rank $r$.

\subsection{Riccati equation} \label{sec:riccati}

The Riccati differential equation is given by
\begin{equation}\label{eq:riccati}
  \dt{X}(t) = A^T X(t) + X(t) A + C^T C - X(t) S X(t), \quad X(0) = X_0,
\end{equation}
where $X \in \R^{m \times m}$, $A \in \R^{m \times m}$, $C \in \R^{k \times m}$, and $S \in \R^{m \times m}$. We note that this is no longer an ODE with an affine vector field and hence our theoretical results do not apply here. As already studied in~\cite{ostermann_convergence_2019}, we take $S = I$ and $A$ is the spatial discretization of the diffusion operator
$$\mathcal D = \partial_x (\alpha(x) \partial_x(\cdot)) - \lambda I$$
on the spatial domain $\Omega = [0,1]$. Furthermore, we take $ \alpha(x) = 2 + 2 \cos(2 \pi x) $ and $\lambda = 1$.
The discretization is done by the finite volume method, as described in~\cite{gander_numerical_2018}. The tall matrix $C \in \R^{k \times m}$ is obtained from $k$ independent vectors $\{1, e_1, \ldots, e_{(k-1)/2}, f_1, \ldots, f_{(k-1)/2} \}$, where
\begin{equation}
  e_i(x) = \sqrt{2} \cos(2 \pi k x) \quad \text{and} \quad f_i(x) = \sqrt{2} \sin(2 \pi k x), \quad i=1, \ldots, (q-1)/2,
\end{equation}
are evaluated at the grid points $\{x_j \}_{j=1}^m$ with $x_j = \frac{j}{m+1}$. The time interval is $[0, T] = [0, 0.1]$.

\begin{figure}
  \centering
  \begin{subfigure}{0.32\textwidth}
    \includegraphics[width=\textwidth]{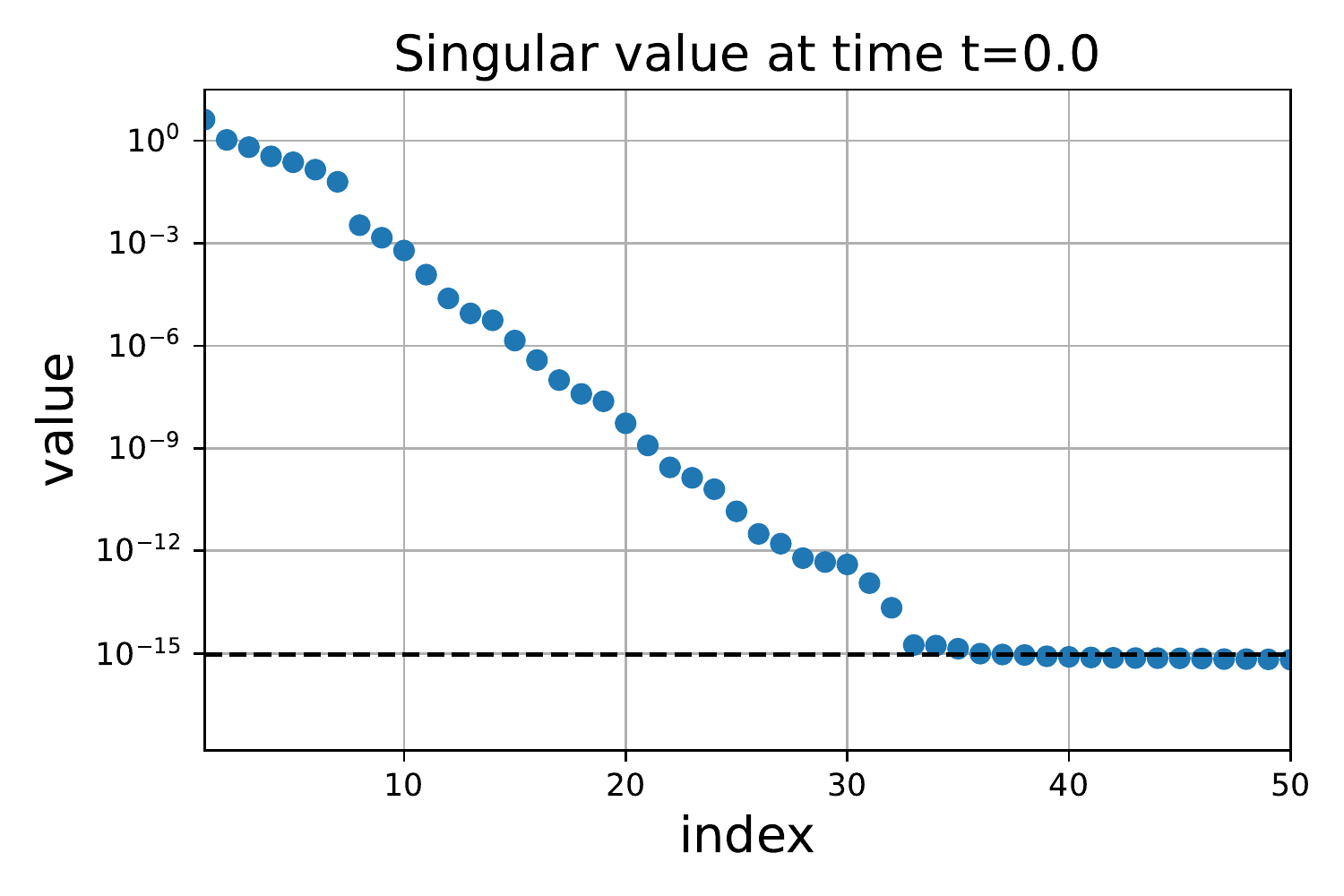}
  \end{subfigure}
  \begin{subfigure}{0.32\textwidth}
    \includegraphics[width=\textwidth]{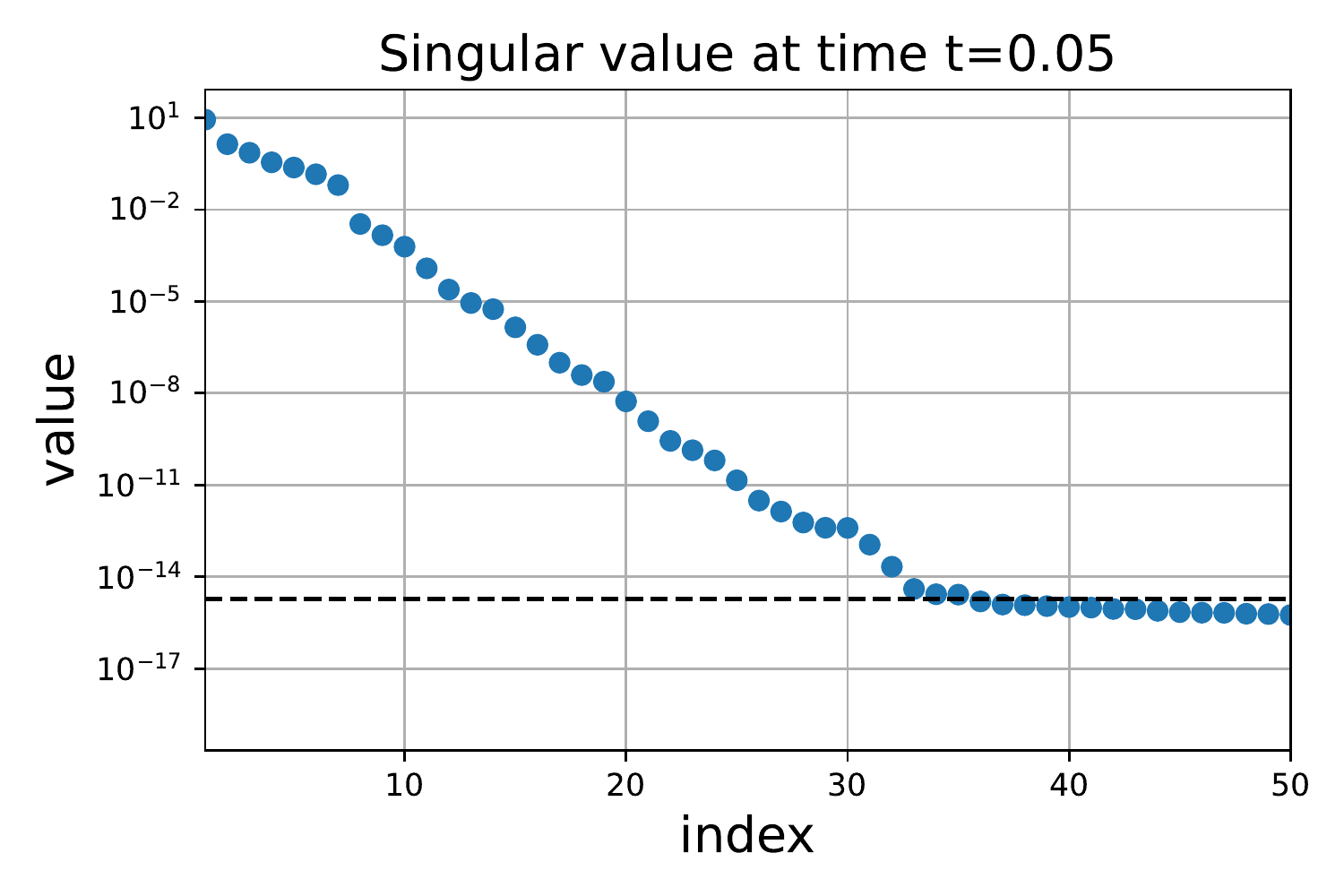}
  \end{subfigure}
  \begin{subfigure}{0.32\textwidth}
    \includegraphics[width=\textwidth]{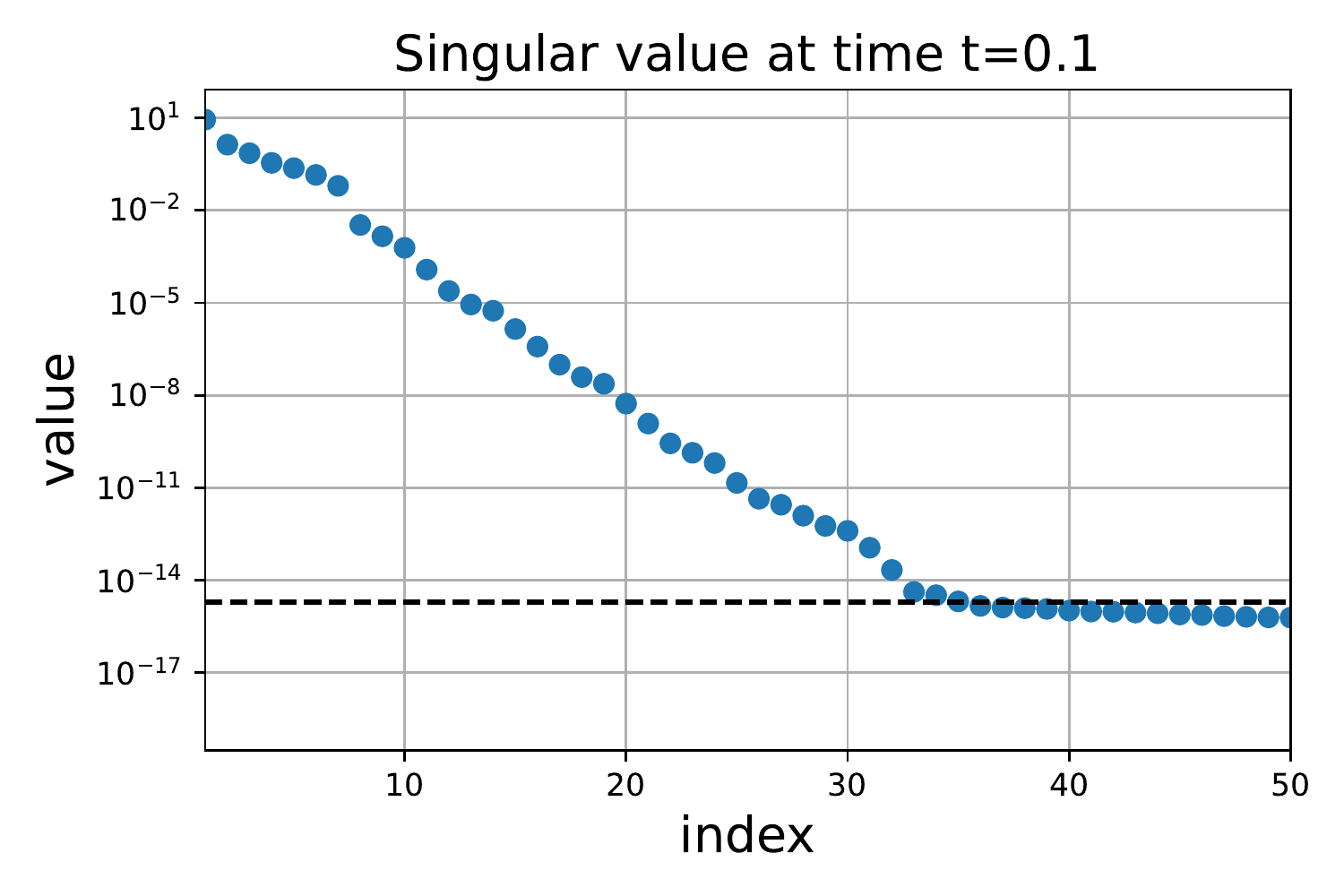}
  \end{subfigure}
  \caption{Singular values of the solution over time of the Riccati ODE~\eqref{eq:riccati}.}
  \label{fig:riccati_solution}
\end{figure}

\begin{figure}
  \centering
  \begin{subfigure}{0.49\textwidth}
    \includegraphics[width=\textwidth]{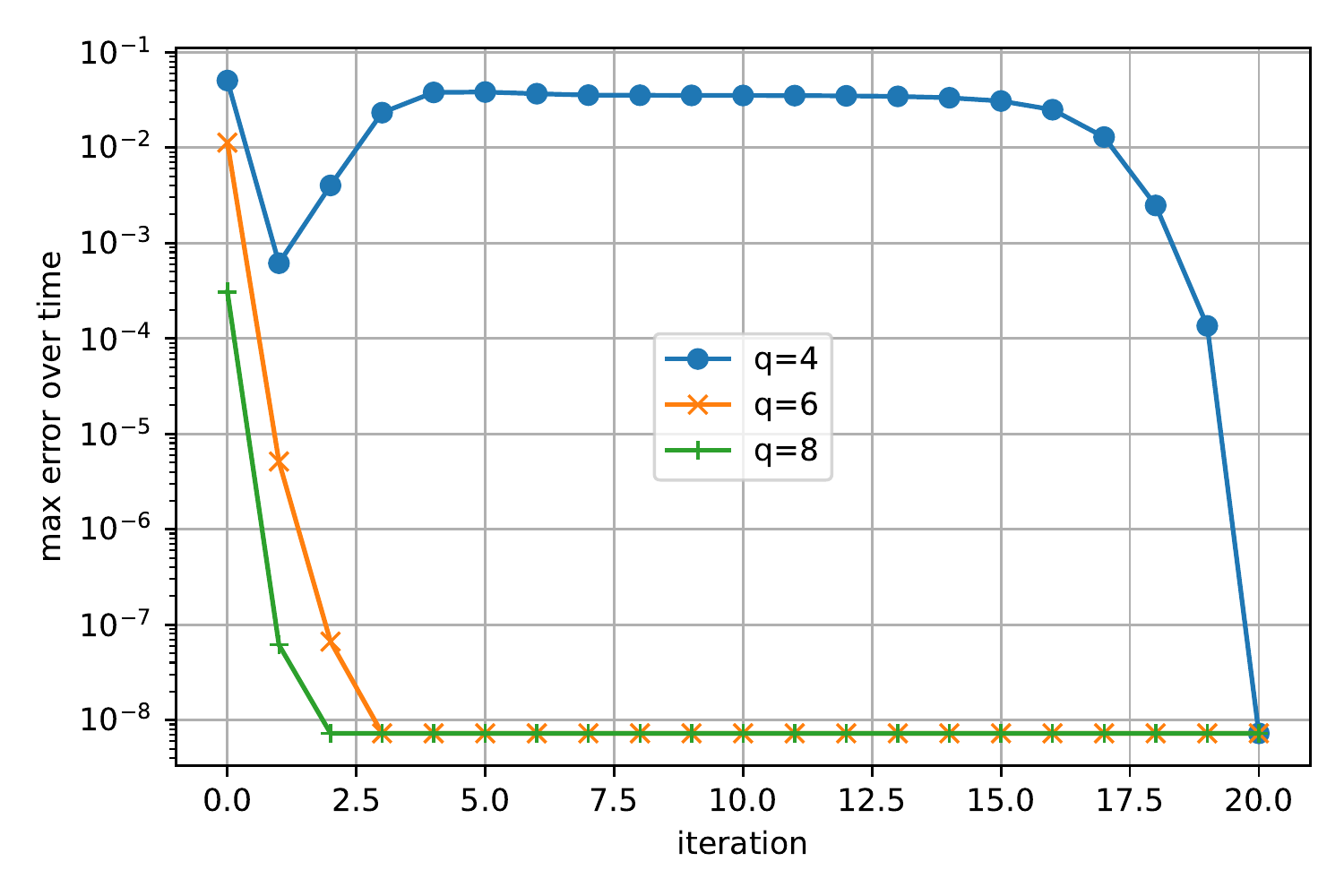}
    \caption{Several coarse ranks with fine rank $r=18$.}
  \end{subfigure}
  \begin{subfigure}{0.49\textwidth}
    \includegraphics[width=\textwidth]{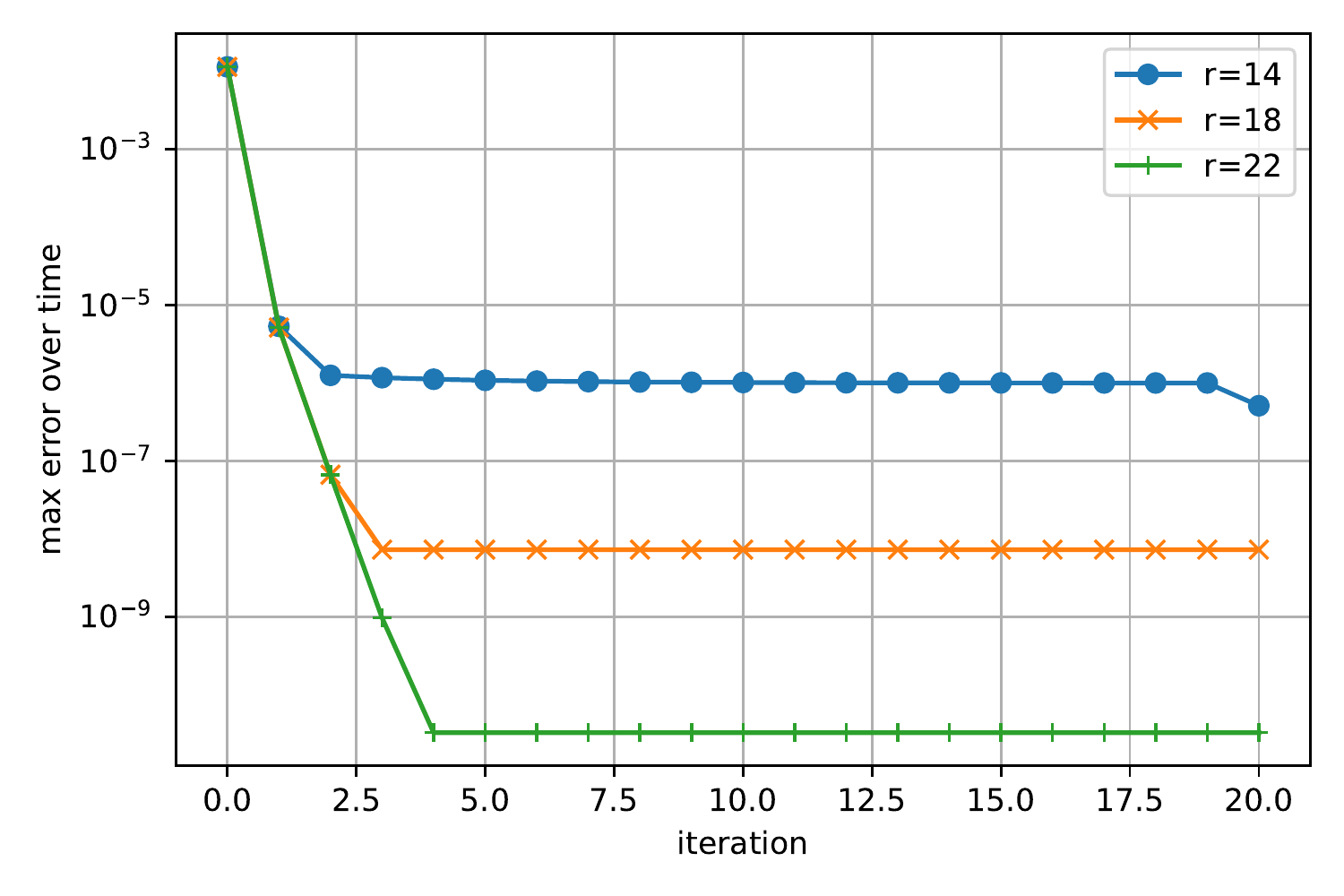}
    \caption{Several fine ranks with coarse rank $q=6$.}
  \end{subfigure}
  \caption{Convergence of the error of low-rank Parareal for the Riccati problem~\eqref{eq:riccati}. Influence of the coarse and fine ranks.}
  \label{fig:riccati_several_ranks}
\end{figure}

As for the other problems, the singular values of the solution (shown in Figure~\ref{fig:riccati_solution}) indicate that we can expect good low-rank approximations on $[0,T]$. We choose the fine rank $r=18$. The convergence of low-rank Parareal is shown in Figure~\ref{fig:riccati_several_ranks}. Unlike the previous problems, the coarse rank $q$ has a more pronounced influence on the behavior of the convergence. While our theoretical results do not hold for this nonlinear problem, we still see that low-rank Parareal converges linearly when the coarse rank $q$ is sufficiently large ($q=6$, $q=8$). The convergence is slower (but still superlinear) when $q=4$. This could be due to the non-constant gaps in the singular values. The influence of the fine rank $r$ is more like for the linear problems.

\subsection{Rank-adaptive algorithm}

Since the approximation rank of the solution is usually not known a priori, it is more convenient for the user to supply an approximation tolerance than an approximation rank. Even though the rank can change to satisfy the tolerance during the truncation steps, Algorithm \ref{def:low_rank_parareal} can be easily reformulated for such a rank adaptive setting. The key idea is to fix the coarse rank to keep the cost of the coarse solver low, while the fine rank is determined by a fine tolerance. 

\begin{definition}[Adaptive low-rank Parareal] \label{def:adaptive low rank parareal}
  Consider a small fixed rank $q$ and a fine tolerance $\tau$. The adaptive low-rank Parareal algorithm iterates
  \begin{align}
    & \text{(Initial value)} \quad          &  & Y_0^k = Y_0, \label{def:adaptive initial value}                                                                                      \\
    & \text{(Initial approximation)}  \quad &  & Y_{n+1}^0 = \flowq \circ \Tq (Y_n^0) + \mathcal E_n, \label{def:adaptive initial approximation}                                                              \\
    & \text{(Iteration)} \quad              &  & Y_{n+1}^{k+1} = \psi^h_{\rank(\mathcal T_{\tau}(Y_n^k))} \circ \mathcal T_{\tau}(Y_n^k) + \flowq \circ \Tq (Y_n^{k+1}) - \flowq \circ \Tq (Y_n^k), \label{def:adaptive iteration}
 \end{align}
 where the notation is similar to that of the previous Def.~\ref{def:low_rank_parareal}, except for $\mathcal T_{\tau}$ which represents the rank-adaptive truncation. In particular,  $\mathcal T_{\tau}(Y)$ is the best rank $q$ approximation of $Y$ so that the $(q+1)$th singular value of $Y$ equals the tolerance $\tau$. The matrices $\mathcal E_n$ are small perturbations, randomly generated such that $\rank(Y_{n+1}^0) = \rank(Y_0)$ and its smallest singular value is larger than the fine tolerance $\tau$. 
\end{definition}

\begin{figure}
  \centering
  \begin{subfigure}{0.49\textwidth}
    \includegraphics[width=\textwidth]{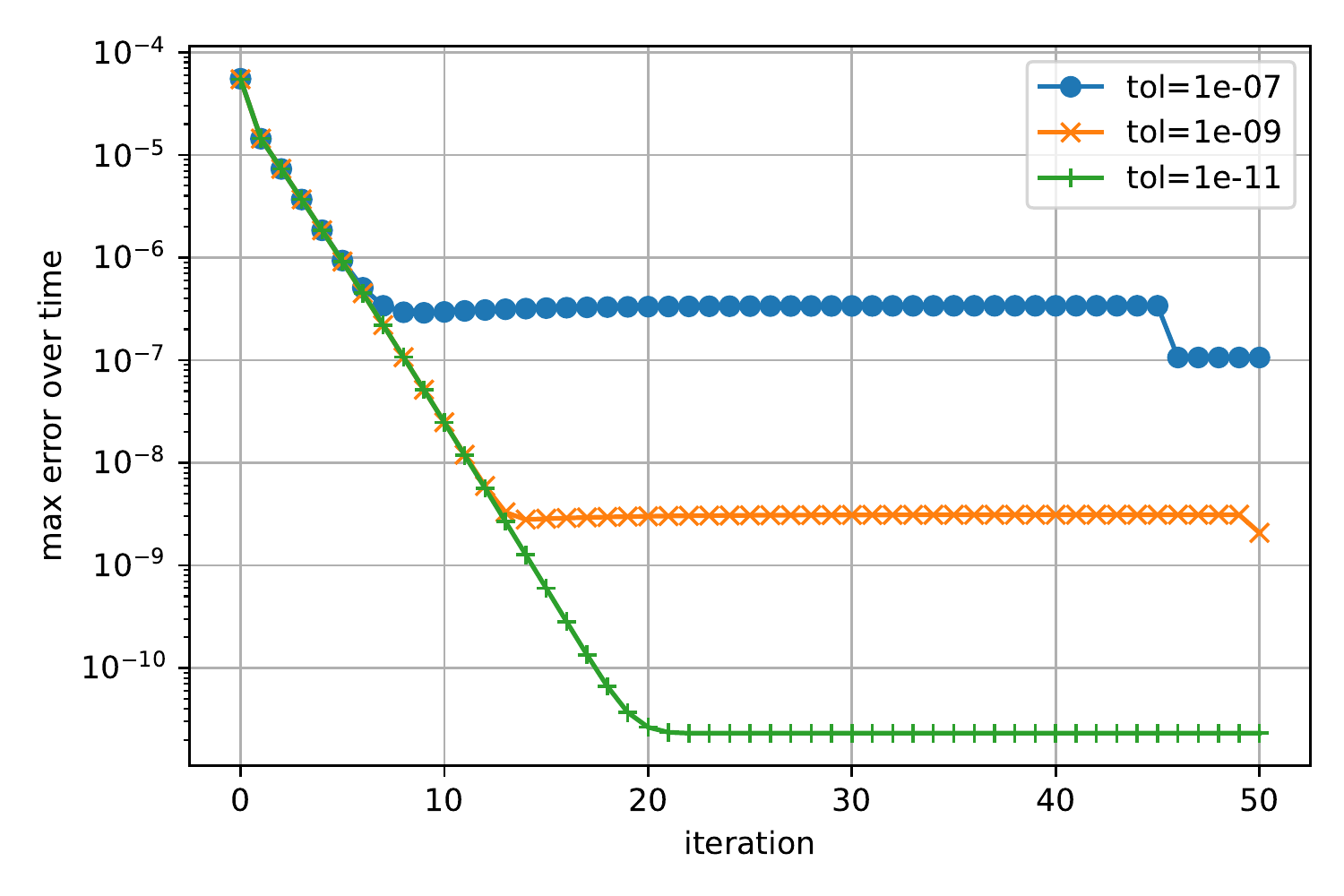}
    \caption{Several tolerances with coarse rank $q=4$.}
    \label{fig: several tolerances}
  \end{subfigure}
  \begin{subfigure}{0.49\textwidth}
    \includegraphics[width=\textwidth]{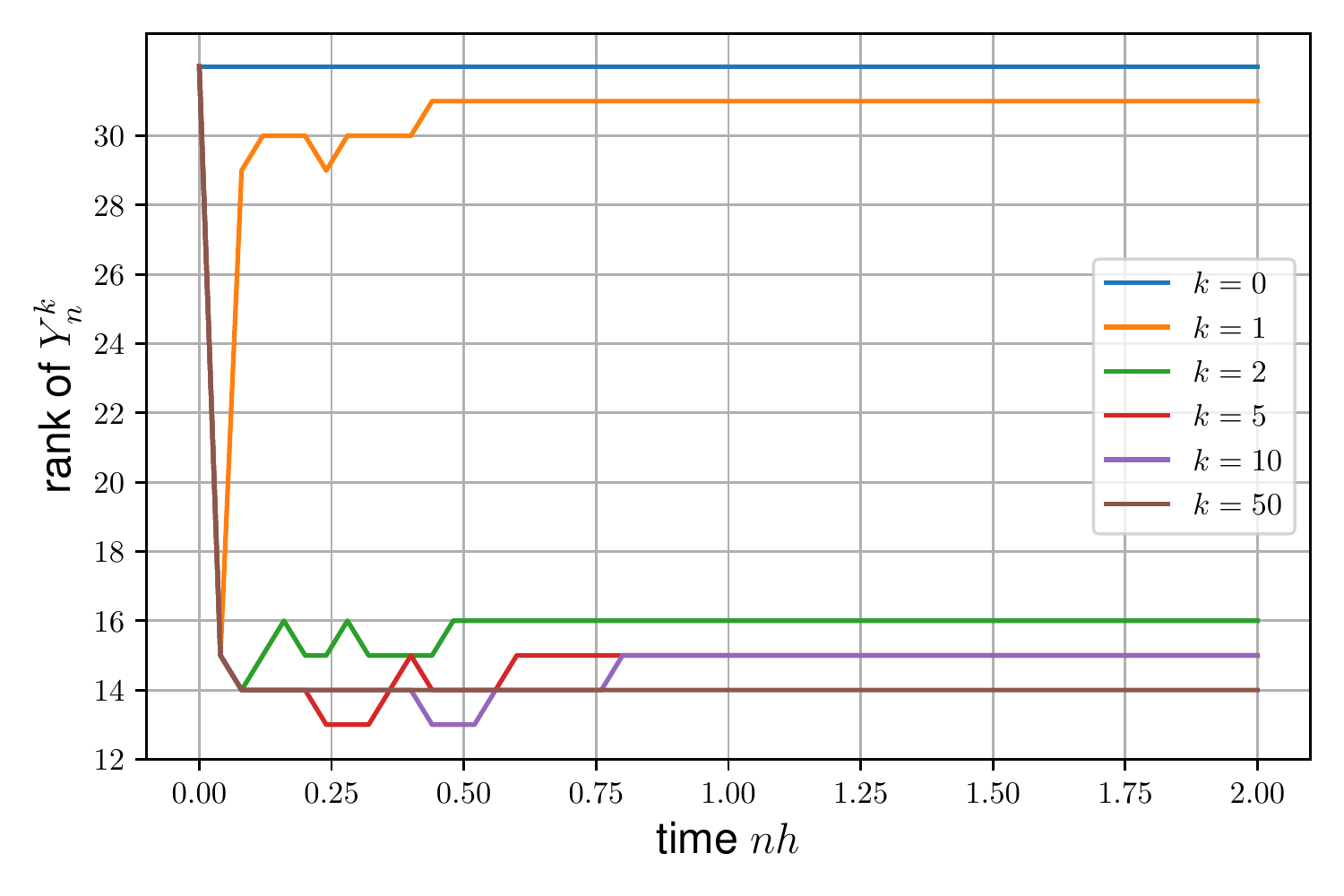}
    \caption{Coarse rank $q=4$ and tolerance =$1e-9$.}
    \label{fig: rank over time and iteration}
  \end{subfigure}
  \caption{Adaptive low-rank Parareal. On the left, the algorithm is applied with several tolerances. On the right, the rank of the solution over time is shown for several iterations.}
  \label{fig: adaptive low-rank parareal}
\end{figure}

Figure \ref*{fig: adaptive low-rank parareal} shows the numerical behavior of this rank-adaptive algorithm. As we can see, the algorithm behaves as desired. Figure \ref{fig: several tolerances} shows the algorithm applied with several tolerances and is comparable to Figure \ref{fig:lyapunov_several_fine_ranks} with several fine ranks. Figure \ref{fig: rank over time and iteration} shows the rank of the solution over time. Already after two iterations, the rank is reduced to almost the numerical rank of the exact solution and the rank does not change much for the rest of the iterations.

\section{Conclusion and future work} \label{sec:conclusion}

We proposed the first parallel-in-time algorithm for integrating a dynamical low-rank approximation (DLRA) of a matrix evolution equation. The algorithm follows the traditional Parareal scheme but it uses DLRA with a low rank as coarse integrator, whereas the fine integrator is DLRA with a higher rank. Taking into account the modeling error of DLRA, we presented an analysis of the algorithm and showed linear convergence as well as superlinear convergence under common assumptions and for affine linear vector fields, up to the modeling error.

In our numerical experiments, the algorithm behaved well on diffusive problems, which is similar to the original Parareal algorithm. Due to the significant difference in computational cost for the fine and coarse integrators, it is reasonable to expect good speed-up in actual parallel implementations. A proper parallel implementation to verify this claim is a natural future work. It may however be more appropriate to first generalize more efficient parallel-in-time algorithms, like Schwarz waveform relaxation and multigrid methods~\cite{gander_time_2018}, to DLRA.

Since DLRA can also be used to obtain low-rank tensor approximations~\cite{lubich_time_2017}, another future work is to extend low-rank Parareal to tensor DLRA. Finally, our theoretical analysis assumes that the ODE has an affine vector field. Since this assumption was only needed in one step of the proof of Lemma~\ref{lem:iteration_error}, it might be possible that it can be relaxed to include certain non-linear vector fields.

\bibliographystyle{abbrvnat}
\bibliography{refs_arxiv2}

\end{document}